\documentclass[10pt]{article} 
 
\usepackage{amsmath} 
\usepackage{amsthm} 
\usepackage{amssymb} 
\usepackage{amsfonts} 
\usepackage{epsfig}
 
\input xypic

\newtheorem{The}{Theorem}[section] 
 
\newtheorem{Pro}[The]{Proposition}

\newtheorem{Lem}[The]{Lemma} 
\newtheorem{Def}[The]{Definition}

\def\endproof{\relax\ifmmode\expandafter\endproofmath\else 
\unskip\nobreak\hfil\penalty50\hskip.75em\hbox{}\nobreak\hfil\bull 
{\parfillskip=0pt \finalhyphendemerits=0 \bigbreak}\fi} 
\def\endproofmath$${\eqno\bull$$\bigbreak} 
\def\bull{\vbox{\hrule\hbox{\vrule\kern3pt\vbox{\kern6pt}\kern3pt\vrule}\hrule}} 
 
\def\cancel#1#2{\ooalign{$\hfil#1\mkern1mu/\hfil$\crcr$#1#2$}}

\def\dirac{\mathpalette\cancel\partial} 
\def\Dirac{\mathpalette\cancel D}

\newcommand{\C}{{\mathbb C}} 
\newcommand{\R}{{\mathbb R}} 
\newcommand{\Z}{{\mathbb Z}} 
\newcommand{\Q}{{\mathbb Q}} 
\renewcommand{\P}{{\mathbb P}}
 
\newcommand{\A}{{\cal A}}

\newcommand{\M}{{\cal M}}

\newcommand{\la}{\langle} 
\newcommand{\ra}{\rangle} 
\newcommand{\ba}{\begin{eqnarray}} 
\newcommand{\na}{\end{eqnarray}} 
\newcommand{\beq}{\begin{equation}} 
\newcommand{\eeq}{\end{equation}}  
\newcommand{\s}{\mathfrak{s}} 
 
\newcommand{\spinc}{\mathrm{Spin}^c}

\title{Exact triangles in monopole homology and the Casson-Walker invariant}  
\author{Matilde Marcolli and Bai-Ling Wang} 
\date{}
 
\begin{document} 
\maketitle 
 
\tableofcontents
 
\section{Introduction} 

The purpose of this paper is to give a general outline of the problem
of the exact triangles in Seiberg--Witten--Floer theory. We present
here the most general case, where the problem consists of producing a
surgery formula relating the monopole homology of a compact oriented
3--manifold $Y$ with an embedded knot $K$, and the monopole homologies
of some 3--manifolds obtained by Dehn surgery on $K$. 

In the series of papers \cite{CMW} \cite{MW2} \cite{MW3} \cite{MW4},
we studied the problem in the case of an integral homology 
3-sphere $Y$, and the 3--manifolds $Y_1$  
and $Y_0$ obtained by Dehn surgery on $K$ with framing $1$ and $0$,
respectively. 

The results of \cite{CMW} \cite{MW2} \cite{MW3} \cite{MW4} are, at
this stage, still to be considered as ``work in progress'', where some
of the proofs need more rigorous presentations. The main
result of that series of papers is that the
Seiberg-Witten-Floer homologies of $Y$ and of the manifolds $Y_1$   
and $Y_0$ are related by an exact triangle
{\small
\[
\diagram
HF^{SW}_{*} (Y_1, g_1) \rto^{w^1_{*}}
 & HF^{SW}_{*}(Y, g_0, \mu)
\dlto_{w^0_{*}} \\
 \bigoplus_{k\in \Z} HF^{SW}_{(*)}(Y_0, \s_k) \uto_{w_{(*)}}
& \\
\enddiagram
\] }
In this triangle the maps $w_*^1, w_*^0$ and $w_{(*)}$ are induced 
by the surgery cobordisms connecting $Y_1$ and $Y$,
$Y$ and $Y_0$, $Y_0$ and $Y_1$, respectively, and $\mu$ is the surgery
perturbation that simulates the effect of Dehn surgery. 

In this paper, we explain how, following the same strategy for the
proof of the surgery formula which we have introduced in the previous
work, we may be able to extend this exact triangle to the more
general case of any closed oriented 3-manifold with a smoothly
embedded knot. 

In the last section of this paper we give a topological application of
the kind of arguments that lead to the proof of the ``geometric
triangles'', namely the surgery formula for monopole homology viewed
at the level of generators. We show that a modified
version of the Seiberg-Witten invariant agrees with the Casson-Walker
invariant, for any closed and oriented rational homology 3-sphere.

Let Y be a  closed oriented 3-manifold, 
with a smooth embedded knot $K$, let $\nu(K)$
be the  tubular neighbourhood of $K$ in $Y$.
Choose an identification of $\nu(K)$ with 
$D^2\times S^1$:
\ba
\nu(K) \cong D^2\times S^1,
\label{frame}
\na 
where $K$ is mapped to the core of the
solid torus $D^2\times S^1$. Under the identification (\ref{frame}), 
on the boundary $T^2$, we fix a basis $m, l$ of $H_1(T^2, \Z)$ such that
$l$ is the longitude (parallel to $K$ under the identification (\ref{frame})) 
and $m$ is the right-handed meridian (intersecting  $l$ once), the orientation
determined by $m\wedge l$ coincides with
the orientation induced from $Y$. The corresponding
longitude and meridian in the knot complement
$ V = Y\backslash \nu (K) $
are denoted by $m', l'$ respectively. 
 Similarly, let $m''$ and $l''$ be the meridian
and longitude in the tubular neighbourhood of the knot $\nu(K)$. The
meridian $m''$ bounds a disk $D^2$ in $\nu(K)$, and $l''$ generates
$H_1(\nu(K),\Z)$ and parallels to $K$. Let $p$ and $q$ be two relatively prime
integers, the Dehn surgery with coefficient $p/q\in \Q\cup \{\infty\}$
on $K$ is the operation of removing $\nu(K)$ and gluing in $D^2\times S^1$ 
by an orientation reversing diffeomorphism $f_{p/q}$
of $T^2$ that satisfies
\[
f_{p/q} (m'') = pm'-ql'.
\]
The resulting manifold is denoted by $Y_{p/q}$. Note that in general 
$Y_{p/q}$ depends on the choice of the identification (\ref{frame}).

Let $\s$ be  a $\spinc$ structure  on $Y$. We shall see that, with a
suitable choice of metrics and perturbation, $(Y, \s)$ has
non-empty monopole moduli space only if $\s|_{\nu(K)}$ 
has trivial determinant, hence we shall always assume that
the   $\spinc$ structure $\s$ is trivial around $K$.
If $K$ represents a trivial homology class 
in $H_1(Y, \Z)$, then there is only one $\spinc$ structure  on 
$Y$ which agrees with $\s$ over $Y-\nu(K)$ and $\nu(K)$. 
Suppose that $K$ represents a torsion element of order $n$ in 
$H_1(Y, \Z)$, which means,
\ba
\label{order:n}
\displaystyle{
\frac {|\hbox{Torsion}( H_1(Y, \Z))|}{|\hbox{Torsion}( H_1(Y-\nu(K),  \Z))|}}
=n.
\na
In other words, $n$ is the minimal positive integer such 
that $n[K]$ is homologous to zero in $H_1(Y, \Z)$. 
Then there exists a $\Z_n$-family of $\spinc$ structures 
$
\s\otimes L_k (k=1, \cdots, n)
$
which agree with $\s$ over $Y-\nu(K)$ and $\nu(K)$, where 
$L_k$ is a complex line bundle whose Euler
class is given by $kPD([K])$. If $[K] \neq 0$ in $H_1(Y, \Q)$,
then  there exists a $\Z$-family of $\spinc$ structures 
$
\s\otimes L_k (k\in \Z)
$
which agree with $\s$ over $Y-\nu(K)$ and $\nu(K)$, where 
$L_k$ is a complex line bundle whose Euler
class is given by $kPD([K])$.

Let $Y_1$ be the $(+1)$-surgery on $K$, 
and $Y_0$ be the 0-surgery on $K$, 
we can consider separately the following cases.

\noindent $\bullet$ 
First we assume that $Y$ is a rational homology 3-sphere with
a smoothly embedded knot $K$ representing
a torsion element of order $n$ in $H_1(Y, \Z)$ in the sense of
(\ref{order:n}). Then $Y_1$ is a rational homology
3-sphere, and $Y_0$ is a rational homology $S^1\times S^2$.
Let $\s$ be a $\spinc$ structure on $Y$, which is
trivial on $\nu(K)$. Gluing
the $\spinc$ structures $s|_{Y-\nu(K)}$ and $\s|_{\nu(K)}$
along $T^2$ via different gauge transformations
on $T^2$ results in
a $\Z_n$-family of $\spinc$ structures on $Y$ and $Y_1$,
and a $\Z$-family of $\spinc$ structures on $Y_0$. 
Without confusion, thinking $K$ as the core of 
the attaching solid torus $\nu(K)$,  we denote these structures
on $Y, Y_1$ and $Y_0$ all by
$ \s\otimes L_k (k\in Z),$
where $L_k$ is a complex line bundle of Euler
class $kPD([K])$.\, With this notation, it is understood
that for $Y$ and $Y_1$, and with $n[K]=0$, there is only a
$\Z_n$-family of $\spinc$ structures. 
Then, with a careful choice of metrics and perturbations
on $Y, Y_1$ and $Y_0$ as in Part I \cite{CMW}, we will
obtain the following exact triangles for the Seiberg-Witten-Floer
homologies: 
{\small
\[
\diagram 
 HF^{SW}_{*} (Y_1, \s\otimes L_m, g_1) \rto^{w^1_{*}}  
 & HF^{SW}_{*}(Y,\s\otimes L_m, g, \mu) 
\dlto_{w^0_{*}} \\  
 \bigoplus_{k\in \Z} HF^{SW}_{(*)}(Y_0, \s\otimes L_{nk+p}) \uto_{w_*}  
& \\ 
\enddiagram 
\] }
which holds for any fixed choice of $m$ and $p$ in $\{ 0,\ldots, n-1
\}$, and for a corresponding choice of perturbation, we have 
{\small
\[
\diagram 
\bigoplus_{k=1}^n HF^{SW}_{*} (Y_1, \s\otimes L_k, g_1) \rto^{w^1_{*}}  
 & \bigoplus_{k=1}^n HF^{SW}_{*}(Y,\s\otimes L_k, g, \mu) 
\dlto_{w^0_{*}} \\  
 \bigoplus_{k\in \Z} HF^{SW}_{(*)}(Y_0, \s\otimes L_k) \uto_{w_*}  
& \\ 
\enddiagram 
\] }
In both cases the homomorphisms $w^1_{*}, w^0_{*}$ and 
$w_*$ are induced by the surgery cobordisms.
These exact triangles generalize the results of 
Part I-IV \cite{CMW} \cite{MW2}\cite{MW3} \cite{MW4},
in the sense that the above exact triangle  reduces to the
exact triangle for an integral homology 3-sphere
when $n=1$.

\noindent $\bullet$ 
Now we assume that $(Y, \s)$ is a closed oriented 
3-manifold of $b_1(Y)>0$ and with a $\spinc$ structure $\s$.
 Let $K$ be a knot representing
a torsion homology class in $H_1(Y, \Z)$ of order $n$
in the sense of (\ref{order:n}).
If $\s$ has  non-trivial determinant in the sense
that $c_1(\s) \neq 0$ in $H^2(Y, \Q)$, then
the Seiberg-Witten-Floer homology
for $(Y, \s)$ is $\Z_{2\ell}$-graded, where
$2\ell$ is the multiplicity of $c_1(\s)$ in
$H^2(Y, \Z)/\hbox{Torsion}$, i.e,
$c_1(\s) (H_2(Y, \Z)/\hbox{Torsion}) = 2\ell$. The gluing of the
$\spinc$ structures corresponding to these three different
surgeries gives rise to a $\Z_n$-family of $\spinc$
structures on $Y, Y_1$ and a $\Z$-family of
$\spinc$ structures on $Y_0$, which
we denote by 
$ \s\otimes L_k \qquad (k\in \Z),$
with the convention as before. Then, for any
$\spinc$ structure  $\s\otimes L_k$ ($k=1, \cdots, n)$ on $Y$
and $Y_1$, the corresponding Seiberg-Witten-Floer
homologies are all $\Z_{2\ell}$-graded, while 
for any $\spinc$ structure  $\s\otimes L_k$ ($k\in \Z$) on $Y_0$,
the Seiberg-Witten-Floer
homology $HF^{SW}_* (Y_0, \s\otimes L_k)$ is 
$\Z_{\ell_{[k]}}$-graded, where $\ell_{[k]}$ is the greatest common
factor in $2\ell$ and $2k$. Similar to Part IV \cite{MW4},
the $\Z_{\ell_{[k]}}$-graded homology $HF^{SW}_* (Y_0, \s\otimes L_k)$
can be lifted to a $\Z_{\ell}$-graded homology 
$HF^{SW}_{(*)} (Y_0, \s\otimes L_k)$. For this case, we derive in this 
paper the following exact triangles:
{\small
\ba
\diagram 
 HF^{SW}_{*} (Y_1, \s\otimes L_m, g_1) \rto^{w^1_{*}}  
& HF^{SW}_{*}(Y, \s\otimes L_m, g, \mu) 
\dlto_{w^0_{*}} \\  
 \bigoplus_{k\in \Z} HF^{SW}_{(*)}(Y_0, \s\otimes L_{nk+p}) \uto_{w_{(*)}}.  
& \\ 
\enddiagram 
\label{b_1>0:n:non-torsion:mp}
\na }
which holds for any choice of $m$ and $p$ in $\{ 0,\ldots, n-1 \}$,
and for a corresponding choice of perturbation, we have 
{\small
\ba
\diagram 
\bigoplus_{k=1}^n HF^{SW}_{*} (Y_1, \s\otimes L_k, g_1) \rto^{w^1_{*}}  
& \bigoplus_{k=1}^n HF^{SW}_{*}(Y, \s\otimes L_k, g, \mu) 
\dlto_{w^0_{*}} \\  
 \bigoplus_{k\in \Z} HF^{SW}_{(*)}(Y_0, \s\otimes L_k) \uto_{w_{(*)}}.  
& \\ 
\enddiagram 
\label{b_1>0:n:non-torsion}
\na }

\noindent $\bullet$ 
If $\s$ is a torsion $\spinc$ structure on $Y$, then for any
$k=1, \cdots, n$, 
$HF^{SW}_{*}(Y, \s\otimes L_k, \Z[[t]])$
and $HF^{SW}_{*}(Y_1, \s\otimes L_k, \Z[[t]])$ are $\Z$-graded
with $\Z[[t]]$-coefficient. The reduced
versions $HF^{SW}_{*}(Y, \s\otimes L_k)$
and $HF^{SW}_{*}(Y_1, \s\otimes L_k)$ are obtained by setting
$t=0$. The Seiberg-Witten-Floer homology 
$HF^{SW}_*(Y_0, \s\otimes L_k)$ is $\Z_{2k}$-graded, and
can be lifted to a $\Z$-graded version, denoted by
$HF^{SW}_{(*)}(Y_0, \s\otimes L_k)$. Then the
exact triangles take the same form as (\ref{b_1>0:n:non-torsion:mp})
and (\ref{b_1>0:n:non-torsion}).

\noindent $\bullet$ 
The remaining case is when a smoothly embedded
knot $K$ represents a non-trivial 
homology class in $H_1(Y, \Q)$. Let $n$ be the 
minimal positive integer such that there is
a  2-cycle in $H_2(Y, \Z)$  intersecting
$K$ at $n$ points.  
Then the $0$-surgery on $Y$ along $K$ yields $Y_0$
satisfying $b_1(Y_0) = b_1(Y) -1$. More precisely,
$H_1(Y_0, \Z)$ is obtained
by replacing the $\Z$-component $\Z \la [K]\ra$ of $H_1(Y, \Z)$
by $\Z_n \la [m'']= [l']\ra$. Notice that 
$Y$ can be thought of as the  manifold obtained by
$0$-surgery on $m'\subset Y-\nu(K) \subset Y_0$,
and $Y_1$ can be thought as the result
of $(+1)$-surgery on
$m'\subset Y-\nu(K) \subset Y_0$. Since $m'$ represents
a torsion element of order $n$ in $H_1(Y_0, \Z)$ in the
sense of (\ref{order:n}), we have 
the exact triangles for $(Y, \s, K)$ obtained from
the corresponding exact triangles for $(Y_0, \s, m')$
in the form of (\ref{b_1>0:n:non-torsion:mp})
and (\ref{b_1>0:n:non-torsion}). Thus,
it is enough to establish the exact triangle 
for a general closed oriented 3-manifold 
$Y$ with a smoothly embedded knot representing 
a torsion element of order $n$ in $H_1(Y, \Z)$.

We now summarize the main theorem of this paper.

\begin{The}\label{main:theorem}
Let $(Y, \s)$ be a closed oriented 3-manifold  with
a $\spinc$ structure which is trivial around a
smoothly embedded knot $K$. Assume that $K$ represents
a torsion element of order $n$ in $H_1(Y, \Z)$ in the sense
of (\ref{order:n}). Assume that the
canonical framing of $K$ is given by
an identification of $D^2\times S^1$ with the 
tubular neighbourhood $\nu(K)$ such that $K$ is given
by $\{0\} \times S^1$. Here the parallel simple curve
on $T^2$ provides  the longitude of $K$ and the right handed meridian
is given by $\partial (D^2) \times \{ pt\}$. The orientation
determined by $m\wedge l$ coincides with the orientation induced  from $Y$.
Let $Y_1$ and $Y_0$ be the manifolds obtained, respectively,
by $(+1)$ and $0$ surgery along $K$ in $Y$.
With a careful choice of metrics and perturbations, we obtain the 
following exact triangle induced
by the surgery cobordisms after possible grading shifts:
{\small
\[
\diagram 
\bigoplus_{k=1}^n HF^{SW}_{*} (Y_1, \s\otimes L_k) \rto^{w^1_{*}}  
 & \bigoplus_{k=1}^nHF^{SW}_{*}(Y, \s\otimes L_k) 
\dlto_{w^0_{*}} \\  
 \bigoplus_{k\in \Z} HF^{SW}_{(*)}(Y_0, \s\otimes L_k) \uto_{w_{(*)}}  
& \\ 
\enddiagram   
\] }
Here $\s\otimes L_k$ is the $\spinc$ structure obtained
by tensoring a $\spinc$ structure $\s$ with a complex line
bundle $L_k$ of Euler class $kPD([K])$. Moreover, for any
fixed $m, p\in \{0, \cdots, n-1\}$, we have the
following more refined version of the exact triangle:
{\small
\[
\diagram
HF^{SW}_{*} (Y_1, \s\otimes L_m) \rto^{w^1_{*}}
 & HF^{SW}_{*}(Y, \s\otimes L_m)
\dlto_{w^0_{*}} \\
 \bigoplus_{k\in \Z} HF^{SW}_{(*)}(Y_0, \s\otimes L_{nk+p}) \uto_{w_{(*)}}
& \\
\enddiagram
\] }
Here again the maps are induced by the surgery cobordisms, possibly
after a shift in the grading.
\end{The}

The major technical steps required in the proof are described in
\cite{CMW}\cite{MW2}\cite{MW3}\cite{MW4}.  Thus, in this paper,
we shall address only those issues that are relevant to 
this general exact triangle, while we refer to the previous papers for the
general setting and results. 

In the last section of the paper, we show that a suitably modified
version of the Seiberg--Witten invariant of a rational homology
3-sphere agrees with the Casson--Walker invariant.
For any rational homology 3-sphere $(Y, \s, g)$ with
a $\spinc$ structure $\s$ and a Riemannian metric, 
the counting of the irreducible Seiberg-Witten
monopoles defines the Seiberg-Witten invariant
\[
SW_Y(\s, g) = \# \bigl(\M_Y^*(\s, g)\bigr),
\]
where each irreducible  monopole in $\M_Y^*(\s, g)$ has
a natural orientation from the linearization of the Seiberg-Witten
equations. As studied in \cite{MW1}, $SW_Y(\s, g)$ depends on the metric and
perturbation used in the definition, in order to obtain a topological
invariant, we can modify $SW_Y(\s, g)$ by a metric and perturbation
dependent correction term as follows. Choose any four manifold $X$
with boundary $Y$, such  that $X$ is endowed with a 
cylindrical-end metric modeled on $(Y, g_Y)$. Choose
a $\spinc$ structure $\s_X$ on $X$ which agrees with
$\s$ on $Y$ over the end, and choose a connection $A$
on $(X, \s_X)$ which extends the unique reducible
$\theta_{\s}$ on $(Y, \s)$. Then we set 
\ba
\xi_Y(\s, g) = Ind_\C (\Dirac^X_A) -\displaystyle{\frac{1}{8} }
\bigl( c_1(\s_X)^2 -\sigma (X)\bigr),
\label{correction:term}
\na
where $Ind_\C (\Dirac^X_A)$ is the complex index of the
Dirac operator on $(X, \s_X)$ twisted with
the extending $\spinc$ connection $A$ and $\sigma (X)$ is
the signature of $X$. By the Atiyah-Patodi-Singer index theorem,
$\xi_Y(\s, g)$ is independent of the choice of
$(X, \s_X)$ and $A$, actually, $\xi_Y(\s, g)$ can be expressed as a
combination of the Atiyah-Patodi-Singer eta invariants for the
Dirac operator and signature operator on $(Y, \s)$:
\[
\xi_Y(\s, g) = \displaystyle{-
\frac 14 \eta^{\dirac_{\theta_{\s}}}_Y (0) -  
\frac 18 \eta^{sign}_Y(0).}
\]

The modified version of the  Seiberg-Witten invariant is defined
as
\ba
\hat{SW}_Y(\s) = SW_Y (\s, g) - \xi_Y(\s, g).
\label{hat:SW}
\na
Then we prove the following equivalence between 
$\hat{SW}_Y$ and the Casson-Walker invariant.

\begin{The}
Let $Y$ be a rational homology 3-sphere. Then,
\[
\sum_{\s\in \spinc (Y)} \hat{SW}_Y(\s) = \displaystyle{\frac{1}{2}}
|H_1(Y, \Z)| \lambda (Y),
\]
where $\lambda (Y)$ is the Casson-Walker invariant of $Y$ (cf. \cite{Walker}).
\end{The}

\vskip .2in

{\bf Acknowledgements} BLW likes to acknowledge the paper
of Ozsv\'ath and Szab\'o \cite{OS} on the theta divisor and the Casson-Walker 
invariant which leads to his proof of the equivalence 
of $\hat{SW}_Y$ and the Casson-Walker invariant, hence 
proving the conjecture formulated in \cite{OS} on the  equivalent
between $\hat{SW}_Y$ 
and their $\hat\theta$ invariant for all rational homology 3-sphere. 
BLW is partially supported by Australia Research Council Fellowship.

\section{The geometric triangle}

In this section, we identify the
monopoles on $Y$ with monopoles on $Y_1$ and $Y_0$.  Suppose given a
smoothly embedded knot   
$K$ in $(Y, \s)$, which represents a torsion element of order $n$ in
 $H_1(Y, \Z)$. We can split $Y$ along a torus as in 
\cite{CMW},
\[
Y= V\cup_{T^2} \nu(K).
\]
We choose a metric on $Y$ with a long cylinder
$[-r, r] \times T^2$, and denote the resulting manifold
as
\[
Y(r) = V\cup_{T^2} ([-r, r] \times T^2) \cup_{T^2} \nu(K).
\]
We additionally require that the chosen metric on $Y$ satisfies the
condition of Lemma 
3.18 in \cite{CMW} in the neighbourhood $\nu(K)$ of the knot:
it has non-negative scalar curvature, strictly
positive away from the boundary. This induces a natural
metric on $Y_0(r)$. On $Y_1(r)$, we need to choose a metric which
agrees with the original metric on $Y (r)$ when restricted
to the knot complement $V$. The induced
metric from $Y_1(r)$ in the torus neighbourhood of $\nu(K)$
is the metric described in Lemma 3.21 \cite{CMW}.

With this choice of the metric, the moduli space of monopoles
on $(Y, \s)$ is non-empty only if $\s|_{\nu(K)} $ is trivial.
Gluing the two $\spinc$ structures $\s|_V$ and $\s|_{\nu(K)}$
along $T^2$ by a gauge transformation on $T^2$ gives rise to a
$\Z_n$-family of $\spinc$ structures on $Y$ and $Y_1$, and to
a $\Z$-family of $\spinc$ structures on $Y_0$. The resulting
$\spinc$ structures  can be classified as the result of tensoring the
original $\spinc$ structure $\s$ with complex
line bundles $L_k (k\in \Z)$ whose
Euler class is given by $kPD([K])$.  The gluing
theorem for 3-dimensional monopoles and the 
perturbation $\mu$ on
$\nu(K)$, ``simulating the effect of surgery'',
provide the decomposition
of the moduli space for
\[
\cup_{k\in \{ 1, \cdots, n\}} \M_Y (\s\otimes L_k).
\]

\begin{The}
\label{decomposition}
With the  choice of perturbations and metrics on $Y, Y_1$ and
$Y_0$ described above, we have the following relation between the
critical sets of the Chern-Simons-Dirac functional
on the manifolds $Y, Y_1$ and $Y_0$:
\[\begin{array}{lll}
&&\bigcup_{k\in \{ 1, \cdots, n\}} \M_Y (\s\otimes L_k)\\[2mm]
&=& \bigcup_{k\in \{ 1, \cdots, n\}} \M_{Y_1}  (\s\otimes L_k)
\cup \bigcup_{k\in \Z} \M_{Y_0}  (\s\otimes L_k).
\end{array}
\]
\end{The}
\begin{proof}
First we assume that $Y$ is a rational homology 3-sphere. When
stretching the neck in $Y(r)$, as $r\to \infty$, we get two manifolds,
each endowed with an infinite cylindrical end,
\[
V(\infty) = V\cup_{T^2} ([0, \infty) \times T^2)\]\[
\nu(K) (\infty) = \nu(K) \cup_{T^2} ((-\infty, 0] \times T^2).
\]
The Seiberg-Witten monopole moduli space of $Y(r)$, for sufficiently large
$r$, can be described in terms of the
moduli spaces of $V(\infty)$ and $\nu(K) (\infty)$ as
analyzed in \cite{CMW}. 

Notice that the moduli spaces of flat connections on $T^2$, modulo
the subgroups of the gauge transformations on $T^2$ which
can be extended to the whole manifolds $V$ or $\nu(K)$, define
the following character varieties:
\[
\chi_0(T^2, V)
= H^1(T^2, \R)/2n\Z \la PD([l])\ra,
\]\[
\chi_0(T^2, \nu(K)) 
= H^1(T^2, \R)/2\Z \la PD([m])\ra.
\]
Thus the character variety $\chi_0(T^2, Y)$ is
$ \chi_0(T^2, Y) = H^1(T^2, \R).$
The covering maps between these character varieties are illustrated as 
follows
{\small \[
\diagram
2\Z \langle PD([m]) \rangle \drto &  & 2n\Z \langle PD([l]) \rangle\dlto \\
& \chi_0(T^2,Y) \dlto^{\pi_1}\drto^{\pi_2}  & \\
\chi_0(T^2,V)\drto &  & \chi_0(T^2,\nu(K))\dlto \\
&  \chi(T^2)  & \\
\enddiagram
\]}
where the maps $\pi_i$ are the covering maps with  fibers as indicated.

Based on the analysis of the space $\M_V^* (\s|_V)$ of irreducible
monopoles on a 3-manifold with 
a cylindrical end modelled on $T^2$, as in \cite{CMW},
we see that the asymptotic limit map defines a continuous
map:
\[
\M_V^* (\s|_V) \stackrel{\partial_{\infty}}{\to}
\chi_0(T^2,V). 
\]
The reducibles on $V$ embed in the character variety $\chi_0(T^2,V)$.
Then, for a sufficiently large $r$, the gluing theorem gives a 
diffeomorphism:
\[
\#_Y : {\cal M}^*_V\backslash \partial_\infty^{-1}(U_\theta)
\times_{\chi_0(T^2,Y)}\chi (\nu(K)) 
\to\bigcup_{k=1}^n {\cal M}^*_{Y(r)}(\s\otimes L_k),
\]
here $U_\Theta$ is a small neighbourhood of the ``bad points'' $\Theta$ 
in $\chi_0(T^2, Y)$ where the twisted Dirac operator
on $T^2$ has non-trivial kernel, and $ \chi (\nu(K))$  
are the reducible lines in $\chi_0(T^2,Y)$. The above
fibred product is obtained (cf.\cite{CMW}) by taking the 
pullbacks of the images of the boundary value maps  under the
projections 
{\small
\[
\diagram
& \chi_0(T^2,Y) \dlto^{\pi_1}\drto^{\pi_2}  & \\
\chi_0(T^2,V) &  & \chi_0(T^2,\nu(K)). \\
\enddiagram
\]}

Let $(u, v)$ be the coordinates on
$\chi_0(T^2,Y) \cong H^1(T^2, \R)$ corresponding
to the holonomy around the longitude $l$ and
the meridian $m$ respectively. In the gluing map $\#_Y$ above, 
$\chi(\nu(K))$ corresponds to the lines
$\{v=2k, k=1, \cdots, n\}$. For each line $\{v=2k\}$,
the image of the gluing map gives
a diffeomorphism onto $\M^*_{Y(r)} (\s\otimes L_k)$.

For each $\spinc$ structure $\s\otimes L_k$, there is a
unique reducible monopole 
on $Y(r)$, which is given by the 
intersection of $\chi(V, \s)$, the flat connections
on $(V, \s)$, with the line $\{ v= 2k\}$ in $\chi_0(T^2, Y)$:
\[
\theta_Y(k) = \{ u=u(\s)\} \times_{\chi_0(T^2, Y)} \{ v= 2k\}.
\]
Here $u(\s)$ is the holonomy of the flat connections
in $\chi(V, \s)$ around the longitude $l'$.

Similarly, we have gluing maps for monopoles on $Y_1(r)$ and
$Y_0(r)$, respectively.
In the gluing map $ \#_{Y_1}$
for $ {\cal M}^*_{Y_1(r)}(\s\otimes L_k)$,
\[
 \#_{Y_1} : {\cal M}^*_V \backslash
\partial_\infty^{-1}(U_\theta)\times_{\chi_0(T^2,Y_1)}\chi (\nu(K)) 
\to\bigcup_{k=1}^n {\cal M}^*_{Y_1(r)}(\s\otimes L_k), 
\]
$\chi (\nu(K))$ is identified with $v-u = 2k+1$.
Similarly, we have the gluing map $ \#_{Y_0}$
for $ {\cal M}^*_{Y_0(r)}(\s\otimes L_k)$,
\[
 \#_{Y_0} : {\cal M}^*_V \backslash
\partial_\infty^{-1}(U_\theta)\times_{\chi_0(T^2,Y_0)}\chi (\nu(K)) 
\to \bigcup_{k\in\Z}{\cal M}^*_{Y_0(r)}(\s\otimes L_k), 
\]
where $\chi (\nu(K))$ is given by $ u = 2k$. For each $k\in \{1,
\cdots, n\}$, the reducible monopole for $(Y_1, \s\otimes L_k)$,
consists of the unique point 
\[
\theta_{Y_1}(k) = \{ u=u(\s)\} \times_{\chi_0(T^2, Y)} \{ v=u+ 2k+1\}.
\]

For $Y_0$ with a non-trivial $\spinc$ structure 
$\s\otimes L_k$ ($k\in \Z, k\neq 0$), the set of reducibles is empty
for any generic perturbation, and for 
$\s\otimes L_0= \s$, it consists of one circle of reducibles
$u=0$ in the cylinder $\chi_0(T^2,Y_0) = \chi_0(T^2,V)$, which
can be perturbed away by introducing a small perturbation as in
Theorem 6.13 \cite{CMW}.

Now we apply the perturbation to simulate the effect of Dehn
surgery. This amounts to a careful choice of perturbation as
in Section 6 \cite{CMW}, which we now briefly describe.

Choose a compactly supported 2-form $\mu$  representing the generator
of $ H^2_{cpt} (D^2 \times S^1)$, defined as in Lemma  3.18 \cite{CMW},
such that we have
\ba
\int_{D^2 \times \{pt \}}\mu  = 1 \label{normalization}
\na
 for any point on $S^1$. Under the isomorphism $H^2_{cpt} (\nu(K))
\cong H_1(\nu(K))$, given by Poincar\'e duality, this form corresponds
to the generator $[\mu]= PD_{\nu(K)}(l)$.
The class of $\mu$ in $H^2( D^2 \times S^1)$
is trivial, and we can write $\mu= d\nu$, where $\nu$ is a 1-form
satisfying $\int_{S^1 \times \{pt \}}\nu  = 1$, i.e.
$\nu= PD_{T^2}(l)$.
Choose on $\nu(K)$ a metric as in Lemma 3.18 \cite{CMW}.

Fix a $U(1)$-connection $A_0$ representing the trivial connection on $T^2$,
For any $U(1)$-connection $A$, define  
$T_A$ to be
\[
T_A(z)= -i\int_{\{z\in D^2\} \times S^1} (A-A_0).
\]

For any  given $\epsilon >0$, we can choose a function $\hat f:\R\to\R$ with
the following properties.

(a) $\hat f$ is continuously differentiable on $(-1,1)$ and satisfies the
periodicity $\hat f(t+2)=\hat f(t)$

(b) the derivative $\hat f'$ has range $\hat f'(t)\in [-1,1]$ for all $t\in
[-1,1]$, and satisfies $\hat f'(1-t) = \hat f'(1+t)$ for $t\in\R$.

(c) the following estimate holds:
$\sup_{t\in [-1+\epsilon, 1-\epsilon]} | \hat f'(t)-t | <
\epsilon.$

Now, for the $\spinc$ structure $\s\otimes L_k$ ($k=1, \cdots, n$),
consider the function $f_k'(t)= \hat f'(t+1)+2k$ and define a
perturbation of the Seiberg-Witten equations
on $(Y, \s\otimes L_k)$ in the following
way:
\ba
\left\{\begin{array}{l}
 F_A = *\sigma(\psi,\psi) + f'_k(T_A) \mu \\[2mm]
\dirac_A (\psi ) = 0
\end{array}
\right. .
\label{p SW}
\na

With respect to the chosen metric on $\nu(K)$, with sufficiently large
positive scalar curvature on the support of $\mu$ as specified in
Lemma 3.18 \cite{CMW}, 
the only solutions
of the perturbed monopole equations are reducibles $(A, 0)$, that satisfy
\ba
F_A = f_k'(T_A)\mu.
\label{deformflat}
\na

In addition to this surgery perturbation, we consider another
perturbation of the Seiberg--Witten equations on the
tubular neighbourhood $\nu(K)$ in $Y$, $Y_1$, and $Y_0$. This
perturbation has the effect of producing a global shift in the
character variety to avoid the bad points on $H^1(T^2, \R)$ when
we deform the unperturbed geometric triangles in 
$\chi (T^2)$ to the perturbed geometric triangles in 
$\chi (T^2)$.

Let $\mu$ be a compactly supported 2-form in $D^2\times S^1$
satisfying (\ref{normalization}). Let $\eta >0$ be some small real
parameter. Consider an additional perturbation
\ba F_A= *\sigma(\psi,\psi) \pm\eta \mu \label{eta:pert} \na
of the curvature equation on $\nu(K)$ inside $Y$ and  inside $Y_0$.
 
This perturbation has the effect of shifting the asymptotic values by
an amount $\eta$. We choose the sign so that the line of reducibles on
$\nu(K)\subset Y$ for $\s\otimes L_k$ 
becomes $\{ (u,v) | v=2k +\eta \}$, the line of reducibles on
$\nu(K)\subset Y_1$ for $\s\otimes L_k$
 remains the same $\{ (u,v) | v-u =2k+ 1 \}$, and the line
of reducibles on $\nu(K)\subset Y_0$ for $\s\otimes L_k$
 becomes $\{ (u,v) | u=2k+ \eta \}$.

On $\nu(K)$ inside $Y$ for $\s\otimes L_k$
 we shall consider the perturbed curvature
equation 
\ba  F_A= *\sigma(\psi,\psi) +
(f'_k(T_A-\eta)+\eta)\mu. \label{eta:mu:pert} \na 

Therefore, we can partition
the moduli spaces for $(Y, \s\otimes L_k)$ ($k=1,\cdots, n$) into the union
of the moduli spaces for $(Y_1, \s\otimes L_k)$ ($k=1,\cdots, n$)
and $(Y, \s\otimes L_k)$ ($k\in \Z)$,  as in Theorem 6.3
\cite{CMW}. This completes the proof of the 
theorem for the case of rational homology 3-sphere $Y$ with
a knot $K$ representing a torsion element of order n in $H_1(Y, \Z)$.

For a general 3-manifold $Y$ with a smoothly embedded
knot $K$ representing a torsion element of order n
in  $H_1(Y, \Z)$, the proof is essentially the same as the case
of rational homology spheres discussed above, and we omit the details
here. 

\end{proof}

\begin{figure}[ht]
\epsfig{file= 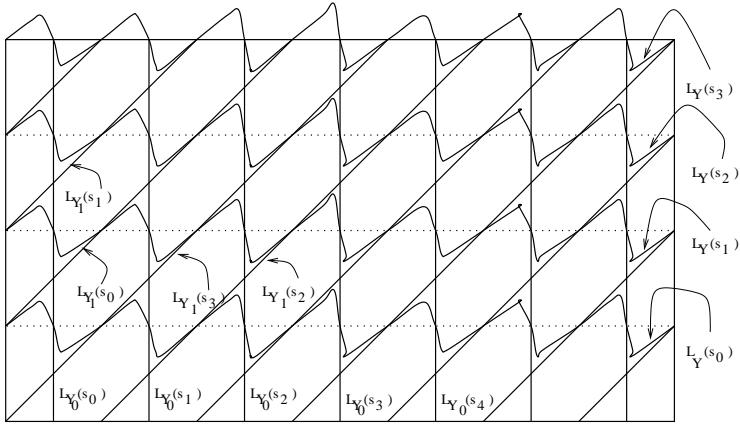,angle= 0}
\caption{The perturbed geometric triangle as in Theorem \ref{decomposition}}
\label{deform:triangle1}
\end{figure}

The perturbation can be illustrated as in Figure
\ref{deform:triangle1} where $n=4$. From now on, we will
use the following notations to denote the reducibles lines
for $Y$, $Y_1$ and $Y_0$ repectively:
\ba\begin{array}{c}
 L_Y(\s_k)= \{ (u,v) | \, v=f'_k(u-\eta) +\eta \}\\[2mm]
 L_{Y_1}(\s_k)=\{ (u,v) | \, v=u+2k+1 \}\\[2mm] 
 L_{Y_0}(\s_k)= \{ (u,v) | \, v=2k+\eta \}. 
\end{array}
\label{reducible:lines}
\na

With a more careful study of the perturbed geometric triangles, we have
the following decomposition of 3-dimensional monopoles
on $Y$ under the Dehn surgery.

\begin{figure}[ht]
\epsfig{file= 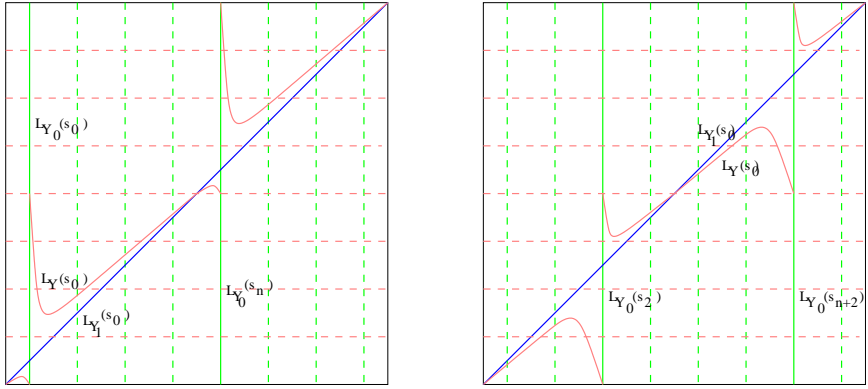,angle= 0}
\caption{The perturbed geometric triangle as in Theorem
\ref{refine:decomposition} 
\label{deform:triangle2}} 
\end{figure}

\begin{The}
With the perturbations and metrics on $Y, Y_1$ and
$Y_0$ as in Theorem \ref{decomposition}, and for any
fixed $m, p \in \{0, \cdots, n-1\}$, there exists
 a further perturbation on $Y$,
such that we have the following relation between the
critical sets of the Chern-Simons-Dirac functional
on the manifolds $Y, Y_1$ and $Y_0$:
\ba \label{split-mp}
 \M_Y (\s\otimes L_m)
= \M_{Y_1}  (\s\otimes L_m)
\cup \bigcup_{k\in \Z} \M_{Y_0}  (\s\otimes L_{nk+p}).
\na
\label{refine:decomposition}
\end{The}
\begin{proof}
In the proof of Theorem \ref{decomposition}, we know
that the perturbed Seiberg-Witten monopoles on $Y, Y_1$
and $Y_0$ are given by the following gluing models (here 
we assume that $Y$ is a rational homology 3-sphere):
\[
{\cal M}^*_{Y(r)}(\s\otimes L_k) 
\cong {\cal M}^*_V\backslash \partial_\infty^{-1}(U_\theta)
\times_{\chi_0(T^2,Y)} \{ v=2k \},
\]
\[
{\cal M}^*_{Y_1(r)}(\s\otimes L_k) 
\cong {\cal M}^*_V\backslash \partial_\infty^{-1}(U_\theta)
\times_{\chi_0(T^2,Y_1)} \{ v= u+2k+1\},
\]
\[
{\cal M}^*_{Y_0(r)}(\s\otimes L_k) 
\cong {\cal M}^*_V\backslash \partial_\infty^{-1}(U_\theta)
\times_{\chi_0(T^2,Y_0)} \{ u=2k\}.
\]
Note that the additional perturbation (\ref{eta:pert}) of the 
 curvature equation on $\nu(K)$ inside $Y$ and  inside $Y_0$
introduce a shift of coordinates $(u, v)$
to $(u+\eta, v+\eta)$. We can introduce these new coordinates,
still denoted by $(u, v)$. Then the reducible line for $V$
is given by $u = -\eta$ in $H^1(T^2, \R)$. 

We will show that there exists a further surgery
perturbation on $Y$ that suits the purpose of
identifying monopoles on $Y, Y_1$ and
$Y_0$ as stated in the Theorem.  
Without loss of the generality,
after possible coordinates change, we can assume that
$\s\otimes L_m = \s\otimes L_0$. 
Fix $p\in \{0, \cdots, n-1\}$. We will construct a function
$f_0': \R \to \R$, 
which depends on a small $\epsilon >0$, such that,
as $\epsilon \to 0$, the curve $v=f_0'(u)$ approaches
the union of lines
\[ L_{Y_1}(\s\otimes L_0)\cup \bigcup_{k\in \Z} L_{Y_0}(\s\otimes
L_{nk+p}) \] 
where
$ L_{Y_1}(\s\otimes L_0)=\{ v=u+1\}, 
 L_{Y_0}(\s\otimes L_{nk+p})= \{ u= 2nk+2p\}.$

We identify $\chi_0(T^2, V)$ with the fundamental domain
\[
\{ u\in \R\} \times \{ 0\le v < 2n \},
\]
in $H^1(T^2, \R)$. The asymptotic values
$\partial_\infty (\M^*_V) \subset \chi_0(T^2, V)$ can
be lifted to $H^1(T^2, \R)$ periodically. Using this
$2n$--periodicity, we only need to construct
a function $f_0: [-1, 2n-1] \to [0, 2n]$,
which depends on $\epsilon$, such that, for any given
$\epsilon\leq \epsilon_0$, we have
\[
\sup_{t\in [-1, 2n-1]\backslash [-\epsilon , \epsilon]}
|f'_0 (t) - t | < \epsilon.
\]
Such function can be easily constructed as in the 
proof of Theorem \ref{decomposition}. Then 
over $[2nk-1, 2nk+2n-1]$, $f_0'(t)$ is defined
to be $f'_0 (t)+2nk$.  See Figure 
  \ref{deform:triangle2} where the perturbation is
illustrated in the cases with $n=4, m=p=0$ and
with $n=4, m=0, p=2$, respectively.  The general case
can be proved by a similar method.
\end{proof}

\section{The relative gradings}

In the previous section, we have the following decomposition:
\[
 \M_{Y, \mu} (\s\otimes L_m)
= \M_{Y_1}  (\s\otimes L_m)
\cup \bigcup_{k\in \Z} \M_{Y_0}  (\s\otimes L_{nk+p})  \]
for any fixed $m, p\in  \{0, \cdots,  n-1\}$. 
Assume that $Y$ is a rational homology 3-sphere.
First we fix a grading on  $\M_{Y, \mu} (\s\otimes L_0)$ defined
in terms of the spectral flow of the linearization operator
for the 3-dimensional Seiberg-Witten equations along a path
connecting an irreducible monopole in $\M_{Y,\mu} (\s\otimes L_0)$
to the unique reducible $\theta_Y(0)$ in the configuration
space for $(Y, \s\otimes L_0)$.  Then
the analysis of the relative grading in Part I section 7 \cite{CMW}
can be applied to induce a compatible grading 
on $\M_{Y_1}  (\s\otimes L_0) \cup 
\bigcup_{k\in \Z} \M_{Y_0}  (\s\otimes L_{nk+p})$ as follows, cf.
Proposition 7.3  -- Corollary 7.7 in \cite{CMW}.

\begin{Pro}\label{relative:grading0} Let $Y$ be a rational homology 3-sphere.
For any fixed $p \in \{0, \cdots, n-1\}$, the Floer complexes
\[
C_*(Y, \s\otimes L_0, \mu) = \oplus _{a\in \M_{Y,\mu} (\s\otimes L_0)}
\Z \la a \ra,
\]
\[
C_*(Y_1, \s\otimes L_0) = \oplus _{a\in \M_{Y_1} (\s\otimes L_0)} \Z
\la a \ra,
\]
\[
C_*(Y_0, \s\otimes L_{nk+p}) = \oplus _{a\in \M_{Y_0} 
(\s\otimes L_{nk+p})} \Z \la a\ra,
\]
have a compatible relative grading of generators in the
following sense.
\begin{enumerate}
\item Suppose given two irreducible critical points $a, b$
in $\M^*_{Y_1}(\s\otimes L_0)$, and the corresponding elements 
$a^\epsilon, b^\epsilon$ in $\M^*_{Y,\mu}(\s\otimes L_0)$
under the above decomposition (\ref{split-mp}). Then
\[
\deg_{Y,\mu}(a^\epsilon)-\deg_{Y,\mu}(b^\epsilon)=
\deg_{Y_1}(a)-\deg_{Y_1}(b).
\]
\item Suppose given two monopoles $a, b$
in $\M^*_{Y_0}(\s\otimes L_{nk+p})$, and the corresponding elements 
$a^\epsilon, b^\epsilon$ in $\M^*_{Y, \mu}(\s\otimes L_0)$
under the above decomposition (\ref{split-mp}). Then 
\[
\deg_{Y_0,\s_k}(a)-
\deg_{Y_0,\s_k}(b) = \deg_{Y,\mu}(a^\epsilon)-
\deg_{Y,\mu}(b^\epsilon) \ \ \hbox{ mod } (2nk+2p).
\]
Therefore, the grading $\deg_{Y,\mu}$ defines a $\Z$-valued lift
of the $\Z_{2nk+2p}$-valued relative index on 
$\M^*_{Y_0}(\s\otimes L_{nk+p})$ under the decomposition
of $\M^*_{Y, \mu}(\s\otimes L_0)$.
\end{enumerate}

\end{Pro}

For a general 3-manifold $(Y, \s)$ with $b_1(Y)> 0$
and  a knot $K$ representing
a torsion element of order $n$, if $\s$ is a torsion
$\spinc$ structure, we know that,
for any $k \in \{0, \cdots, n-1\}$, the $\spinc$ structure 
$\s\otimes L_k$ also has a torsion class $c_1(\s\otimes L_k)$. Thus,
after a small perturbation to get rid of the $(S^1)^{b_1(Y)}$-family
of reducibles, we have a $\Z$-graded
$$\M^*_{Y, \mu}(\s\otimes L_k) \cong \M_{Y,\mu}(\s\otimes L_k).$$
Then it is easy to see that the results
in Part I section 7 \cite{CMW} hold in this case as well 
without any substantial change.

Now assume that $c_1(\s)$ is a non-torsion element, with multiplicity $2\ell$
in $H^2(Y, \Z)/\hbox{Torsion}$. Then, for
any $k \in \{0, \cdots, n-1\}$, the set of generators
$$\M^*_{Y, \mu}(\s\otimes L_k) \cong \M_{Y, \mu}(\s\otimes L_k)$$
is $2\ell$-graded, and so is 
$$\M^*_{Y_1}(\s\otimes L_k) \cong \M_{Y_1}(\s\otimes L_k).$$
Then for any $k\in \Z$, any non-empty moduli space
$$\M^*_{Y,\mu}(\s\otimes L_k) \cong \M_{Y, \mu}(\s\otimes L_k)$$
is $\Z_{2\ell_{[k]}}$-graded, where $2\ell_{[k]}$
is the maximum common factor of $2\ell$ and $2k$.
We can choose a relative grading on $\M^*_{Y, \mu}(\s\otimes L_0)$
by the spectral flow of the linearization operator along
a path connecting any irreducible monopole
in $\M^*_{Y}(\s\otimes L_0)$ to a fixed monopole
$a_0$ in $\M^*_{Y, \mu}(\s\otimes L_0)$.  
Note that this grading is $2\ell$-graded. Then the analysis 
in Part I section 7 can also be applied to obtain the
following proposition.

\begin{Pro}\label{relative:grading1}
For $p \in \{0, \cdots, n-1\}$ and $k\in \Z$, the Floer complexes
$C_*(Y, \s\otimes L_0), C_*(Y_1, \s\otimes L_0)$, and 
$C_*(Y_0, \s\otimes L_{nk+p})$ have a compatible 
relative grading of generators in the following sense.
\begin{enumerate}
\item Suppose given two irreducible critical points $a, b$
in $\M^*_{Y_1}(\s\otimes L_0)$, and the corresponding elements
$a^\epsilon, b^\epsilon$ in $\M^*_{Y, \mu}(\s\otimes L_0)$, 
then
\[
\deg_{Y,\mu}(a^\epsilon)-\deg_{Y,\mu}(b^\epsilon)=
\deg_{Y_1}(a)-\deg_{Y_1}(b),
\]
as a $\Z_{\ell}$-valued function.
\item Suppose given two monopoles $a, b$
in $\M^*_{Y_0}(\s\otimes L_{nk+p})$, and the corresponding elements
$a^\epsilon, b^\epsilon$ in $\M^*_{Y, \mu}(\s\otimes L_0)$,
then
\[
\deg_{Y_0,\s_k}(a)-
\deg_{Y_0,\s_k}(b) = \deg_{Y,\mu}(a^\epsilon)-
\deg_{Y,\mu}(b^\epsilon) \ \ \hbox{ mod } (2\ell_{[nk+p]}).
 \]
Thus, the induced grading $ \deg_{Y,\mu}$ on 
$\M^*_{Y_0}(\s\otimes L_{nk+p})$, under the
decomposition (\ref{split-mp}), defines a choice of a
$\Z_{2\ell}$-valued lifting of the $\Z_{2\ell_{[nk+p]}}$-valued
relative index. 
\end{enumerate}
\end{Pro}

In the rest of this section we discuss the induced
gradings on the various $\M_Y^*(\s\otimes L_m)$, for $m \in \{1,
\cdots, n-1\}$. 
Notice that, in the case of a rational homology 3-sphere $Y$,
the perturbed line $L_{Y}(\s\otimes L_m)$ is a
deformation of the line $L_{Y}(\s\otimes L_0)$. Along this
deformation, the corresponding parameterized
reducibles form a path
connecting $\theta_Y(0)$ to $\theta_Y(m)$. As the following Lemma
shows, this deformation
can be realized as a perturbation of the monopole equations
for $(Y, \s\otimes L_0)$. 

\begin{Lem} 
For any $m \in \{1, \cdots, n-1\}$, there exists a perturbation
$\mu_m$ for the Seiberg-Witten monopole equations
on $(Y, \s\otimes L_0)$,
 such that there is a diffeomorphism:
\[
\M_{Y, \mu_m} (\s\otimes L_0) \cong \M_Y(\s\otimes L_m).
\]
Moreover, if $Y$ is a rational homology 3-sphere, 
the grading on $\M_Y^*(\s\otimes L_m)$
induces a grading on $\M_{Y, \mu_m}^* (\s\otimes L_0)$. Under the
above identification,
the resulting grading differs from the original grading
on  $\M_{Y}^* (\s\otimes L_0)$ by
the wall-crossing formulae studied in
\cite{MW1} for $(Y, \s\otimes L_0)$. If $b_1(Y)> 0$,
then the induced grading on $\M_Y(\s\otimes L_m)$
from $\M_{Y, \mu_m} (\s\otimes L_0)$ agrees with the
relative grading on $\M_Y(\s\otimes L_m)$. 
\label{pert-0m}
\end{Lem}

\begin{proof}
The first claim follows from the gluing models
of the monopoles in $\M_{Y} (\s\otimes L_0)$
and $ \M_Y (\s\otimes L_m)$. Then, using the results in \cite{MW1}, we
know that, in the 
case of the rational homology 3-sphere,
the spectral flow of the twisted Dirac operator
along the path of reducibles gives the 
index shift on
$\M_{Y}^* (\s\otimes L_0)$ and $\M_{Y, \mu_m}^* (\s\otimes L_0)$
according to the wall-crossing formulae derived in \cite{MW1}. Again
by the results of \cite{MW1},
for $Y$ with $b_1(Y)>0$, the induced
grading from the parametrized spectral flow is same as the
original relative index on $ \M_Y (\s\otimes L_m)$. 
\end{proof}

Thus, Lemma \ref{pert-0m} provides a consistent way of assigning a
choice of absolute grading on the various $\M_Y (\s\otimes L_m)$.
The degree shift of Lemma \ref{pert-0m} can be described as follows. 
Let $\theta_Y(t)$ be the path of reducibles 
on $(Y, \s\otimes L_0) $ for the family of perturbations
connecting $\M_{Y}^* (\s\otimes L_0)$
to $\M_{Y, \mu_m}^* (\s\otimes L_0)$,
then the wall crossing formulae in \cite{MW1}
tell us that
the index shift is given by the complex spectral flow
\ba
SF_\C(\dirac_{\theta_Y(t)}). 
\label{index:shift}
\na

\section{Geometric limits and the holomorphic triangles}

\subsection{Surgery cobordisms}

We first briefly describe the surgery cobordisms from
$Y_1$ to $Y$, from $Y$ to $Y_0$, and from $Y_0$ to $Y_1$,
respectively, as in \cite{MW3}.  The cobordism $W_1$, from
$Y_1$ to $Y$, is obtained by removing from the trivial cobordism 
$Y_1\times [0, 1]$ an $S^1\times D\cong \nu(K)\times \{ 1\}$, where $D$ is
a disk, and $\nu(K)$ is the tubular neighbourhood of the knot in
$Y_1$, and then attaching a 2-handle with framing $-1$.
We denote by $D_1$ the core disk of the 2-handle in $W_1$. 
Similarly, the cobordism $W_0$, form $Y$ to $Y_0$, is obtained by removing
from the trivial cobordism $Y_0\times [0, 1]$ an $S^1\times D\cong
\nu(K)\times \{ 0\}$ and attaching a 2-handle with framing zero.
We denote by $D_0$ the core disk of the 2-handle in $W_0$.
Attaching the two-handle has the effect of modifying the
boundary component $Y_1\times \{ 1\}$ in the trivial cobordism to the
boundary component $Y\times \{ 1 \}$ in the non-trivial cobordism
$W_1$, or, respectively,  the
boundary component $Y_0\times \{ 0 \}$ in the trivial cobordism 
to the boundary component $Y\times \{ 0 \}$ in $W_0$. 
The cobordism $\bar W_2$ connecting $Y_0$ and $Y_1$,
satisfies the relation
$$ \bar W = \bar W_2 \# \C \P^2, $$
where $\bar W$ is the composite cobordism $\bar W=\bar W_0 \cup_Y \bar
W_1$.

We assume that the 3-manifolds $Y_1$, $Y$, and $Y_0$ are endowed with metrics
with a long cylinder $T^2\times [-r,r]$, as specified in
section 2 (see also \cite{CMW}). We consider the manifolds $W_1$ and
$W_0$ endowed 
with infinite cylindrical ends $Y_1 \times (-\infty,
-T_0]$ and $Y\times [T_0,\infty)$, and $Y_0 \times
[T_0,\infty)$ and $Y\times (-\infty, -T_0]$, respectively. 
As in \cite{MW3}, we can decompose the cobordisms $W_i$ as
\ba 
W_i= V\times \R \cup_{T^2\times \R} T^2\times [-r,r]\times \R
\cup_{T^2\times \R}  W_i(\nu(K)). 
\label{split:cobord} 
\na
The  region $W_i(\nu(K))$ has the following property. There 
is a compact set ${\cal K}$ in $W_i$ such that the intersection
${\cal K} \cap W_i(\nu(K))$ is obtained by attaching a 2-handle 
$D\times D$ to the product $\nu(K)\times [-T_0, T_0]$, and, outside
of ${\cal K}$, the region ${\cal K}^c \cap W_i(\nu(K))$ consists of product
regions $\nu(K)\times [ T_0 ,\infty)$ and $\nu(K)\times (-\infty, -T_0]$,
and $T^2\times [r_0,r] \times [-T_0,T_0]$.

As in \cite{MW3},  consider an interior point $x_i$
contained in the core disk of the 2-handle, $x_i \in D_i$,
and we denote by $\hat W_i$ the 
punctured cobordism $\hat W_i = W_i \backslash \{ x_i \}$. Similarly,
we can consider the punctured manifold
$$ \hat W_i(\nu(K)) = W_i(\nu(K))\backslash \{ x_i \}. $$
In  the manifolds $\hat W_i(\nu(K))$, endowed with an extra
asymptotic end of the form $S^3 \times [0,\infty)$ at
the puncture, we can identify a product region 
\ba
 {\cal V}=\nu(K)_{r_0}\times \R \cong D\times (D_i\backslash \{ x_i
\}). 
\label{prodW} 
\na 
Thus, we identify the manifold $W_i$ with a
connected sum
$$ W_i = \hat W_i \# Q_i, $$
with a long cylindrical neck $S^3\times [-T(r),T(r)]$, and with $Q_i$ a
4-ball, where $S^3$ is  decomposed as the union of
two solid tori in the standard way,
$S^3=\nu \cup \tilde \nu$, with $\nu\cong \tilde\nu \cong D\times
S^1$. Then the product region ${\cal V}$ of (\ref{prodW}) 
in $W_i$ identifies the standard solid torus $\nu$ in $S^3$ with the
neighbourhood 
$\nu(K)$ of the knot $K$ in $Y$. Similarly, there is a product
region $\tilde {\cal V}$ which identifies the other solid torus
$\tilde\nu$ in $S^3$ with the 
tubular neighbourhood $\nu(K)$ in $Y_i$, after the surgery.
The resulting punctured cobordism can be written as
\ba
 \hat W_i(r)=( V_r\times \R) \cup {\cal V}(r) \cup \tilde {\cal V}(r).
\label{hat:W(r)}
\na

We now impose a choice of metrics and perturbations for
the Seiberg-Witten equations
on the cobordisms as in subsection 2.2 of \cite{MW3}. Then we can adopt the
results of \cite{MW2}\cite{MW3} to understand
the asymptotic limits of finite energy monopoles, under the
splitting of the punctured cobordisms as $r\to\infty$, as in
(\ref{hat:W(r)}).

\subsection{Geometric limits and holomorphic triangles}

We assume that $Y$ and $Y_1$ are endowed with the
$\spinc$ structure $\s\otimes L_0$. The other $\spinc$ structures can
be studied analogously. 

On the surgery cobordism $W_1(r)$, there
is a $\Z$-family of $\spinc$ structures whose
restrictions to the two ends agree with
$\s\otimes L_0$ on $Y$ and $Y_1$ respectively. 
We will study the moduli spaces of monopoles on
$W_1(r)$ with asymptotic values in
$\M_{Y_1}(\s\otimes L_0)$ and 
$\M_{Y, \mu}(\s\otimes L_0)$ at the two ends,
where $\mu$ is the surgery perturbation defined by $f_0'$ as
in the proof of Theorem \ref{refine:decomposition}.
In the case $b_1(Y)>0$, we only consider
the components of minimal energy, as defined in \cite{MW3} \cite{MW4}, 
among all the possible moduli spaces with different $\spinc$
structures and with the given asymptotic values.
In particular, for a rational homology 3-sphere $Y$, 
we only consider the moduli spaces of minimal
dimension among the $\Z$-family of $\spinc$ structures.

With this convention understood, we denote by
$\M^{W_1}(a_1, a)$ 
the moduli space with asymptotic values
 $a_1\in\M_{Y_1}(\s\otimes L_0)$ and $a\in\M_{Y,\mu}(\s\otimes L_0)$.
Similarly, we denote by $\M^{W_0}(a, a_0)$ 
the moduli space with asymptotic values
 $a\in\M_{Y,\mu}(\s\otimes L_0)$ and 
$a_0 \in \M_{Y_0}(\otimes L_{nk+p})$ for $k\in \Z$, for a
fixed $p\in\{0, \ldots, n-1\}$. Under  a generic choice of the perturbation,
all these moduli spaces $\M^{W_1}(a_1,a)$ and
$\M^{W_0}(a,a_0)$ are cut out transversely and of the expected
dimension.

The convergence and gluing arguments developed in \cite{MW1} can be
applied to this case as well, to give the following 
compactifications of $\M^{W_1}(a_1, a)$ and $\M^{W_0}(a, a_0)$. 

\begin{Pro}\label{compactification}
Suppose that $\M^{W_1}(a_1, a)$ is non-empty, then $\M^{W_1}(a_1, a)$
admits a compactification to a manifold
with corners, where the  
codimension $1$ boundary strata consist of
\ba 
\label{compact:W1}
\begin{array}{c}
\bigcup_{c\in \M_{Y, \mu}^*(\s\otimes L_0)}
\M^{W_1}(a_1,c)\times \hat \M_{Y,\mu}(c,a)
\\[2mm]
\cup \bigcup_{c_1 \in \M_{Y_1}^*(\s\otimes L_0)} \hat\M_{Y_1}(a_1,c_1)\times
\M^{W_1}(c_1,a), 
\end{array} 
\na
and with extra components
\ba 
\begin{array}{c} 
\hat \M_{Y_1}(a_1,\theta_1)\times U(1) \times
\M^{W_1}(\theta_1,a) \\[2mm]
\M^{W_1}(a_1,\theta)\times U(1) \times\hat
\M_{Y,\mu}(\theta,a), 
\end{array} 
\label{red:compact:W1} 
\na
when splitting through the reducibles $\theta_1$ and
$\theta$ in $ \M_{Y_1}^*(\s\otimes L_0)$ and
$\M_{Y, \mu}^*(\s\otimes L_0)$ respectively . 
We also have the similar compactification for 
$\M^{W_0} (a, a_0)$.
\end{Pro}

Now we can describe the geometric limits of monopoles
in  $\M^{W_1}(a_1, a)$ and $\M^{W_0} (a, a_0)$
when stretching $r\to \infty$ inside
\[
 \hat W_i(r)= V_r\times \R \cup {\cal V}(r) \cup \tilde {\cal V}(r).
 \]
We only describe the case of $\hat W_1(r)$. With 
similar arguments we have the corresponding 
geometric limits for $\hat W_0(r)$.
The proof of the results stated below on these geometric limits
follows from the same arguments of \cite{MW2}\cite{MW3}.

Consider elements $a_i^{(1)}\in \M_{Y_1}^*(\s\otimes L_0)$
and $a_j(\epsilon) \in \M_{Y,\mu}^*(\s\otimes L_0)$,
which we can write as
$$ a_i^{(1)}=[(A_i^-,\psi_i^-)\#(a_{\infty,i}^-,0)] $$
$$ a_j(\epsilon)= [(A_j^+(\epsilon),\psi_j^+(\epsilon))\#
(a_{\infty,j}^+(\epsilon),0)], $$
with 
$$[A_i^-,\psi_i^-], \, \, [A_j^+(\epsilon),\psi_j^+(\epsilon)]\in \,
\, \M_V^*(\s\otimes L_0|_V)$$
$$ a_{\infty,i}^-\in \, \, L_{Y_1}(\s\otimes L_0), 
 a_{\infty,j}^+(\epsilon)\in \, \, L_{Y,\mu}(\s\otimes L_0).$$

\begin{The} (Proposition 3.1 and Remark 6.2 in \cite{MW3})
Assume that $\M^{W_1(r)}\bigl(a_i^{(1)}, a_j(\epsilon)\bigr)$
is non-empty for all
sufficiently large $r$. Then a
 family of solutions $[\A_1(r), \Psi_1(r)]$ in 
$\M^{W_1(r)}\bigl(a_i^{(1)}, a_j(\epsilon)\bigr)$
defines the following geometric limits on 
$ V_r\times \R \cup {\cal V}(r) \cup \tilde {\cal V}(r)$
as $r\to \infty$.

\noindent 
 (a). A finite energy solution $[\A',\Psi']^\epsilon$ of the perturbed
Seiberg-Witten equations  on $V \times \R$, with a radial limit
$a_\infty(\epsilon)$ in $\partial_\infty (\M_V^*)\subset
\chi_0(T^2,V)$, and with temporal limits $[A,\psi]_1^\epsilon$ and
$[\tilde A, \tilde \psi]_1^\epsilon$ in
$ \partial_\infty^{-1}(a_\infty(\epsilon)) \subset \M_V^*. $

\noindent 
(b). Two paths $[A(t),\psi(t)]_1^\epsilon$ in $\M_V^*$, for
$t\in [-1,0)$ and $t\in (0,1]$, with
\[
\begin{array}{lll}
 [A(-1),\psi(-1)]_1^\epsilon = [A_i^-, \psi_i^-], &\quad&
 \lim_{t\to 0-} [A(t),\psi(t)]_1^\epsilon =
[A,\psi]_1^\epsilon \\[2mm]
 [A(1),\psi(1)]_1^\epsilon = [ A_j^+(\epsilon) ,
\psi_j^+(\epsilon)], &\quad&
 \lim_{t\to 0+} [A(t),\psi(t)]_1^\epsilon = [\tilde A,
\tilde \psi]_1^\epsilon 
\end{array}\]
These paths induce a continuous, piecewise smooth path 
$a_1^\epsilon(t)$ on $\partial_\infty (\M_V^*)$
satisfying $ a_1^\epsilon(t)=\partial_\infty
\bigl([A(t),\psi(t)]_1^\epsilon\bigr)$, with
$$ a_1^\epsilon (-1)=a_{\infty,i}^-  \ \  a_1^\epsilon
(0)=a_\infty(\epsilon) \ \  a_1^\epsilon (1)=
a_{\infty,j}^+(\epsilon). $$ 
As $\epsilon\to 0$, these geometric limits define paths
$[A(t),\psi(t)]$ and $a(t)$  with
$$ a(-1)=a_{\infty,i}^-, \ \  a(0)=a_\infty, \ \ 
a(1)=a_{\infty,j}^+, $$
in $\partial_\infty (\M_V^*)\subset
\chi_0(T^2,V)$,  and $a_\infty=\lim_\epsilon a_\infty(\epsilon)$.

(c) There is  a holomorphic triangle in $H^1(T^2,\R)$ with
vertices at
$$ a^-_{\infty,i}, \vartheta_1, a^+_{\infty,j}(\epsilon) $$
and sides given by parameterized arcs along the lines 
$L_{Y_1}(\s\otimes L_0)$, $L_{Y,\mu}(\s\otimes L_0)$ and 
$ \{ a_1^\epsilon(t) \}\subset \partial_\infty(\M_V^*).$
Here we denote 
$ \vartheta_1 = L_{Y_1}(\s\otimes L_0)\cap L_{Y,\mu}(\s\otimes
L_0).$ 
\label{geom:lim1} 
\end{The}

\begin{figure}[ht]
\epsfig{file=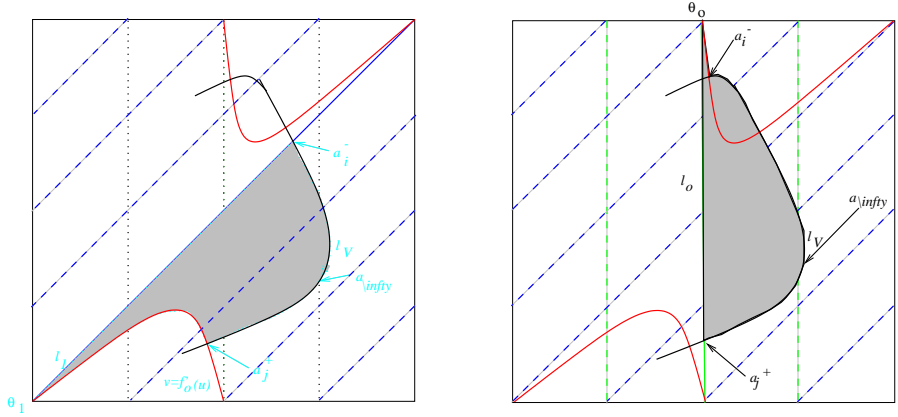,angle= 0}
\caption{Holomorphic triangles for $W_1(r)$ and $W_0(r)$}
\label{holomorphic:triangles}
\end{figure}

This Theorem shows that the 
moduli space $\M^{W_1(r)} \bigl(a_i^{(1)}, a_j(\epsilon)\bigr)$
is characterized by the geometric limits on $V\times \R$ from (a) and
the holomorphic triangles in (c). Two typical
 holomorphic triangles
for $\M^{W_1(r)} \bigl(a_i^{(1)}, a_j(\epsilon)\bigr)$
or  $\M^{W_0(r)} \bigl(a_i(\epsilon), a_j^{(0)}\bigr)$
are illustrated in 
Figure \ref{holomorphic:triangles} where $n=4, m=k=0$ and $p=2$. 
Here the points $\vartheta_i$ are the intersection points
$$ \vartheta_1 = L_{Y_1}(\s\otimes L_0) \cap L_{Y,\mu}(\s\otimes
L_0)$$
for $W_1$ and
$$ \vartheta_0 = L_{Y,\mu}(\s\otimes L_0) \cap L_{Y_0}(\s\otimes
L_{nk+p}) $$
in the case of $W_0$. 
In other words, $\vartheta_i$ is the restriction to
$T^2=\partial 
\nu=\partial \tilde \nu$ of the unique reducible point $\theta_{S^3}$
at the puncture in the cobordism.

In order to glue back the
geometric limits and the holomorphic triangles, there are some
admissible conditions for the geometric limits on $V\times \R$,
which are studied in section 6.1,  section 6.2 and 
section 6.3 \cite{MW3}. 
As in section 6 \cite{MW3}, assuming that the element 
$a_\infty$ is away from the bad
points in the character variety of $T^2$, we denote by
$\hat\M_{V\times \R}(a_\infty)$ the {\em balanced energy} moduli space
 of the Seiberg-Witten equations on $V\times \R$ with 
asymptotic value $a_\infty$ in the radial direction. Then we have
$$ \hat \M_{V\times \R}(a_\infty) = \bigcup_{[A,\psi],[\tilde A,\tilde
\psi]}\hat \M_{V\times \R}( [A,\psi],[\tilde A,\tilde
\psi], a_\infty), $$
with 
$$ [A,\psi],[\tilde A,\tilde\psi] \in \partial_\infty^{-1}(a_\infty)
\subset \M^*_V. $$  
Each moduli space 
$$ \M_{V\times \R}( [A,\psi],[\tilde A,\tilde
\psi], a_\infty), $$
for a fixed choice of $[A,\psi]$ and $[\tilde A,\tilde
\psi]$ in $\partial_\infty^{-1}(a_\infty)$ in $\M_V^*$, is a smooth
finite dimensional oriented manifold of the expected dimension, where
the orientation is given by the corresponding determinant line bundle
of the linearization operator for the monopole equations on 
$V\times \R$.

Remember that we are only considering the components with minimal
energy or minimal dimension among all the moduli spaces
of finite energy monopoles with fixed asymptotic values.
Recall that there is a notion of admissible triples
(cf. Definition 6.4 \cite{MW3}),  which
singles out those elements $([A,\psi],[\tilde A,\tilde \psi], a_\infty)$
which arise as part of the geometric limits of solutions in
$\M^{W_1(r)}(a_1,a)$  
(or $\M^{W_0(r)}(a_1,a_0)$, or $\M_{Y(r)\times \R}(a,b)$ etc).
For example, in the case of $\M^{W_1(r)}(a_1,a)$, with
$$ a_1=[(A^-,\psi^-)\#(a^-,0)] $$
$$ a=[(A^+,\psi^+)\#(a^+,0)], $$ 
 an admissible triple  $([A,\psi],[\tilde A,\tilde \psi], a_\infty)$  
must satisfy the following conditions:
there exists a smooth regular
parameterization $a(t)$, for $t\in [-1,1]$ of the path in
$\partial_\infty(\M_V^*)$ connecting $a^-$ and $a^+$, such that
$a(0)=a_\infty$, and corresponding smooth paths $[A(t),\psi(t)]$ in
$\M_V^*$, for $t\in [-1,0)$ and $t\in (0,1]$, satisfying  
$\partial_\infty[A(t),\psi(t)]=a(t)$, and with
$$ [A(-1),\psi(-1)]=[A^-,\psi^-] \ \ \hbox{ and
} \ \ \lim_{t\to
0_-}[A(t),\psi(t)]=[A,\psi] $$
$$ \lim_{t\to 0_+}[A(t),\psi(t)]=[\tilde A,\tilde \psi] \ \ \hbox{ and
} \ \ [A(1),\psi(1)]=[A^+,\psi^+]. $$
By the results of \cite{MW3}, we know that the possible choices of
$a_\infty$ and of the admissible data in $\hat\M_{V\times
\R}(a_\infty)$ are uniquely determined by the inequivalent holomorphic 
triangles $\Delta$ with vertices
$\{ a^-,\vartheta_1, a^+ \}$ and sides along the union of
Lagrangians $\ell \cup \ell_1 \cup \ell_\mu$, with $\ell$ defined by
the asymptotic values $\partial_\infty(\M_V^*)$. We denote by
$\Xi_{W_1} (a_1, a)$ the set of such inequivalent holomorphic
triangles. 
Moreover, under the identification of the choice of admissible triples 
$([A,\psi],[\tilde A,\tilde \psi], a_\infty)$
with the choice of inequivalent oriented holomorphic triangles in 
$\Xi_{W_1} (a_1, a)$,  the gluing map
gives an orientation preserving diffeomorphism 
\ba \label{glue:W_1}
\#_{W_1}: 
\bigcup_{([A,\psi],[\tilde A,\tilde \psi], a_\infty)\in \Xi_{W_1} (a_1, a)}
\M_{V\times\R}( [A,\psi],[\tilde A,\tilde \psi], a_\infty)
 \to
\M^{W_1}(a_1,a). 
\na 

Similarly, there are orientation preserving diffeomorphisms
given by the gluing maps for $W_0(r)$, $(Y(r)\times \R, \s\otimes L_0)$,
$(Y_1(r)\times \R, \s\otimes L_0)$ and $(Y_0(r)\times \R, \s\otimes
L_{nk+p})$ defined over the set of admissible triples, which are in turn 
determined by the corresponding  inequivalent oriented holomorphic triangles
or holomorphic discs. The corresponding sets
are denoted by 
$\Xi_{W_0}(a, a_0)$, $\Xi_{Y}(a, b)$,  $\Xi_{Y_1}(a_1, b_1)$
and  $\Xi_{Y_0}(a_0, b_0)$, respectively. 
We summarize all the gluing theorems for the
geometric limits and holomorphic triangles (discs) as follows:

\begin{The} \label{final:glue}
Suppose given a pair of monopoles $a^{(1)}, b^{(1)}$ in
$\M_{Y_1}(\s\otimes L_0)$, and a pair of monopoles
$a(\epsilon), b(\epsilon)$ in $\M_{Y_1, \mu}(\s\otimes L_0)$,
where the surgery perturbation $\mu$ determined by $f_0'$ depends on a 
small parameter $\epsilon$. Suppose given a pair of monopoles 
$a^{(0)}, b^{(0)}$ in $\M_{Y_0}(\s\otimes L_{nk+p})$
for some $k\in \Z$ and $p\in \{0, \cdots, n-1\}$. Then, for
sufficiently large $r$,
the gluing maps give the following
orientation preserving diffeomorphisms
\[\begin{array}{l}
\bigcup_{\Delta\in \Xi_{Y} (a(\epsilon), 
b(\epsilon))}
\M_{V\times\R}( \Delta)
\stackrel{\#_{Y(r)\times \R}}{\longrightarrow}  
\M_{Y\times\R} (a(\epsilon), b(\epsilon)),
\\
\bigcup_{\Delta\in \Xi_{Y_1}(a^{(1)}, 
b^{(1)})}
\M_{V\times\R}( \Delta)
 \stackrel{\#_{Y_1(r)\times \R}}{\longrightarrow} \M_{Y_1\times\R} (a^{(1)}, 
b^{(1)}),
\\
\bigcup_{\Delta\in \Xi_{Y_0
}(a^{(0)}, b^{(0)})}\M_{V\times\R}( \Delta)
 \stackrel{\#_{Y_0(r)\times \R}}{\longrightarrow} 
\M_{Y_0\times\R} (a^{(0)}, b^{(0)}),
\\
\bigcup_{\Delta\in \Xi_{W_1} (a^{(1)}, 
a(\epsilon))}
\M_{V\times\R}( \Delta)
\stackrel{\#_{W_1}}{\longrightarrow} \M^{W_1}(a^{(1)}, 
a(\epsilon)),
\\
\bigcup_{\Delta\in \Xi_{W_0} 
(a(\epsilon), a^{(0)})}
\M_{V\times\R}( \Delta)
\stackrel{\#_{W_0}}{\longrightarrow} \M^{W_0} (a(\epsilon), a^{(0)}),
\end{array}\]
where, for simplicity, we denoted a triple $([A,\psi],[\tilde A,\tilde
\psi], a_\infty)$ by $\Delta$. 
\end{The}

\section{Proof of exactness and the surgery triangle}

In this section, we will prove the main theorem of this paper. Notice
that we only prove the refined exact triangle in Theorem
\ref{main:theorem} for $m=0$. The arguments for this case can be
adapted to give an analogous proof for $m\neq 0$.

\subsection{The chain homomorphisms}

Recall that the moduli spaces for $W_1$ and $W_0$ which we consider
are only the components of minimal energy (or dimension). They are
smooth and oriented manifolds of the expected dimension, and they can
be compactified according to Proposition \ref{compactification}. In
particular, whenever one such moduli space is 0-dimensional, we have a 
counting of points with the orientation. Thus,
as in \cite{MW3}, we can define the chain homomorphisms
between the Floer chain groups of $(Y, \mu, \s\otimes L_0)$,
$(Y_1,  \s\otimes L_0)$ and $(Y_0, \s\otimes L_{nk+p})$ as follows.

\begin{Def}
We define the map $w_*^1: C_*(Y_1)\to C_*(Y,\mu)$ by assigning the
matrix elements 
$$ \la a , w_*^1(a_1) \ra =\# \M^{W_1} (a, a_1),$$
where the right hand side is a counting of points with the orientation
if $\M^{W_1} (a, a_1)$ is 0-dimensional and non-empty, and  it is 
$\la a , w_*^1(a_1) \ra =0$ otherwise. 
Similarly, we define the map 
$w_*^0: C_*(Y,\mu)\to \oplus_k C_{(*)}(Y_0, \s\otimes L_{nk+p})$,
with the the matrix coefficients 
$$ \la a_0 , w_*^0(a) \ra =\# \M^{W_0}(a, a_0) $$
if $\M^{W_0}(a, a_0))$ is 0-dimensional and non-empty, and  
$ \la a_0 , w_*^0(a) \ra =0$ otherwise.
Here the chain groups $C_{(*)}(Y_0, \s\otimes L_{nk+p})$ are 
equipped with the lifting of the relative
grading on $\M_{Y_0}(\s\otimes L_{nk+p})$ as described in 
Proposition \ref{relative:grading1}.
\label{w10}
\end{Def} 

Notice that on $W_1$ there is a $\Z$-family of 
$\spinc$ structures which agree with 
$\s\otimes L_0$ when restricted to the two
ends. From Lemma 4.1 in \cite{MW3}, the choice of components with 
minimal energy and dimension essentially 
excludes other $\spinc$ structures when
we consider the moduli space $\M^{W_1}(a_1, a)$ with fixed
asymptotic values $a_1 \in \M_{Y_1}^*(\s\otimes L_0)$ and
$a\in \M^*_{Y, \mu}(\s\otimes L_0)$.

In the case of $b_1(Y)>0$ and $c_1(\s\otimes L_0)$ is non-torsion,
when both  $\M_{Y, \mu}(\s\otimes L_0)$ and 
$\M_{Y_1}(\s\otimes L_0)$ are $\Z_{2\ell}$-graded, with
$2\ell$ the multiplicity of $c_1(\s\otimes L_0)$ in
$H^2(Y, \Z)/\hbox{Torsion}$, then Proposition \ref{relative:grading1}
ensures that there exists
a compatible grading on $\M_{Y_1}(\s\otimes L_0)$
and a  $\Z_{2\ell}$--lifting of the relative grading 
on any $\M_{Y_0}(\s\otimes L_{nk+p})$. For those
$\spinc$ structures $\s\otimes L_0$  with torsion $c_1(\s\otimes L_0)$,
the corresponding Floer homology is defined
to be
\[
HF_*^{SW}(Y, \s\otimes L_0, \Z[[t]])|{t=0},
\]
as described in \cite{MW4}.  Then it is easy to see that
the choice of components with minimal energy and dimension makes 
$w_*^1$ and $w_*^0$ well-defined.  

Then, with the help of the compactifications in
Proposition \ref{compactification},
 the proofs of Lemma 4.3 and Lemma 4.4 in \cite{MW3}
go through without any substantial change to give the following
 Lemma.

\begin{Lem} 
\begin{enumerate}
\item The maps $w_*^i$ are chain homomorphisms.
\item Suppose given $a_1$ in $\M_{Y_1}^*(\s\otimes L_0)$
and $a_0$ in $\M_{Y_0}(\s\otimes L_{nk+p})$,  such that the relative
index induced from $\M_{Y, \mu}^*(\s\otimes L_0)$ as in
Proposition \ref{relative:grading1} is zero, then the 
composite map $w^0_* \circ w^1_*$ is given by
\[
 \la w^0_* \circ w^1_*(a_1), a_0 \ra =
\# \M^{W} (a_1, a_0).
\]
Here we use the same convention on the choice of the
moduli spaces for $W =W_1\#_Y W_0$.
\end{enumerate}
\label{chain:homo}
\end{Lem}

Thus, we have obtained a sequence of
chain complexes induced by the surgery
cobordisms:
\[
 0\to C_*(Y_1, \s\otimes L_0)\stackrel{w^1_*}{\longrightarrow}
C_*(Y,\mu, \s\otimes L_0)\stackrel{w^0_*}{\longrightarrow} 
\oplus_{k\in \Z} C_{(*)}(Y_0,\s\otimes L_{nk+p}) \to 0,
\]
for any fixed $p\in \{0, \cdots, n-1\}$.
We shall prove that this sequence is exact and that the corresponding
exact triangle is the surgery triangle, namely the connecting
homomorphism is also defined via a surgery cobordism. The gluing Theorem
\ref{final:glue}
for the admissible geometric limits on $V\times \R$ 
and the corresponding holomorphic triangles (or discs)
will play a crucial role in the proof of these statements.

\subsection{Proof of exactness}
We first prove  that $w^1_*$ is injective and
$w^0_*$ is surjective, then we show that as in \cite{MW3},
$w^0_* \circ w^1_* =0$. This, together with the specific properties of 
the maps $w^1_*$ and $w^0_*$ will be sufficient to prove the exactness
in the middle term of the sequence:
\[ 0\to C_*(Y_1, \s\otimes L_0)\stackrel{w^1_*}{\longrightarrow}
C_*(Y,\mu, \s\otimes L_0)\stackrel{w^0_*}{\longrightarrow}
\oplus_{k\in \Z} C_{(*)}(Y_0,\s\otimes L_{nk+p}) \to 0. \]

We can partition the moduli space
$\M_{Y, \mu} (\s\otimes L_0)$ according to the
decomposition (\ref{split-mp}) of Theorem \ref{refine:decomposition}, 
\[
\M^*_{Y, \mu} (\s\otimes L_0) \cong 
 \M^*_{Y_1}  (\s\otimes L_0) \cup \bigcup_{k\in \Z} 
\M_{Y_0}  (\s\otimes L_{nk+p}),
\]
where $\mu$ represents the surgery perturbation determined by $f_0'$
as in the proof of Theorem \ref{refine:decomposition}. Here 
$f_0'$ depends on a small parameter $\epsilon >0$. That is,
we identify $\M^*_{Y, \mu} (\s\otimes L_0) $ with 
 a collection of points
\ba
 {\cal M}_{Y,\mu}^* (\s\otimes L_0)
= \{ a^{(1)}_i(\epsilon) \}_{i=1,\ldots r} \cup \{
a^{(0)}_j(\epsilon) \}_{j=r+1,\ldots, s }, 
\label{partition}
\na
so that, as we let $\epsilon \to 0$, the points
$\{ a^{(1)}_i(\epsilon) \}_{i=1,\ldots r} $
get identified with the corresponding elements
$\{ a^{(1)}_i\}_{i=1,\ldots r}$
in $ \M^*_{Y_1}  (\s\otimes L_0)$ and, similarly, the points
$\{ a^{(0)}_j(\epsilon) \}_{j=r+1,\ldots, s }$
get identified with the corresponding elements 
$\{ a^{(0)}_j \}_{j=r+1,\ldots, s }$
in $ \cup_{k\in \Z} \M_{Y_0}  (\s\otimes L_{nk+p})$.

\begin{Lem}\label{inj:surj} The coefficients of the maps $w_*^1$ and
$w_*^0$ satisfy
\begin{enumerate}
\item $ \la a^{(1)}_i(\epsilon), w_*^1(a^{(1)}_j \ra =\delta_{ij}$;
\item $ \la a^{(0)}_i, w_*^0(a^{(0)}_j(\epsilon) \ra
=\delta_{ij}$. \end{enumerate}\noindent
Thus, $w_*^1$ is injective and $w_*^0$
is surjective.
\end{Lem}
\begin{proof} Using the gluing Theorem \ref{final:glue}, we can 
describe the moduli space  
$\M^{W_1}\bigl(a_i^{(1)}, a_j^{(1)}(\epsilon)\bigr)$ for
two monopoles $a_i^{(1)}, a_j^{(1)}(\epsilon)$ of
relative index $0$. With our convention
on the choice of components for these moduli spaces, we know
that $\M^{W_1}\bigl(a_i^{(1)}, a_j^{(1)}(\epsilon)\bigr)$
is zero-dimensional, if non-empty, and obtained
as the gluing of the admissible geometric limits
on $V\times \R$ and the corresponding holomorphic
triangles. 

As we let $\epsilon \to 0$, it is easy to see that the holomorphic 
triangles degenerate to certain holomorphic discs, and
the admissible geometric limits
on $V\times \R$ for $\M^{W_1}\bigl(a_i^{(1)}, a_j^{(1)}(\epsilon)\bigr)$
are identified with the admissible geometric limits
on $V\times \R$  for $\M_{Y_1\times \R}(a_i^{(1)}, a_j^{(1)})$.
Since the relative index of $a_i^{(1)}, a_j^{(1)}$ is
zero, $\M_{Y_1\times \R}(a_i^{(1)}, a_j^{(1)})$,
being a zero-dimensional moduli space of minimal energy,
is empty unless $a_i^{(1)}=a_j^{(1)}$, in which case
$\M_{Y_1\times \R}(a_i^{(1)}, a_i^{(1)})$ consists of a unique
solution. This proves that
\[
\la a^{(1)}_i (\epsilon), w_*^1(a^{(1)}_j\ra =\delta_{ij}.
 \]
Similarly, we obtain
$\la a^{(0)}_i, w_*^0(a^{(1)}_j(\epsilon) \ra =\delta_{ij}.$\end{proof}

Now we prove the exactness in the middle term. We proceed
as in \cite{MW3} to show that $w_*^0\circ w_*^1=0$, which, together
with Lemma \ref{inj:surj} is
sufficient to establish the exact triangle.  Again, we
will use heavily the gluing theorem \ref{final:glue}
to analyze the moduli spaces on the 
cobordisms. The following Lemma is the direct
consequence of the results of Lemma \ref{chain:homo} and Lemma
\ref{inj:surj}.

\begin{Lem}
Suppose given  $a^{(1)}_i\in \M_{Y_1}(\s\otimes L_0)$ and
$a^{(0)}_j(\epsilon) \in \M_{Y, \mu}(\s\otimes L_0)$ which
corresponds to $a^{(0)}_j \in \cup_{k\in \Z} \M_{Y_0}(\s\otimes L_{nk+p})$. 
The coefficients of 
the composition  map $w_*^0\circ w_*^1$ satisfy
\[
\la a^{(0)}_j, w_*^0\circ w_*^1 (a^{(1)}_i) \ra
= \la a^{(0)}_j(\epsilon), w_*^1 (a^{(1)}_i) \ra
+ \la a^{(0)}_j, w_*^0 (a^{(1)}_i(\epsilon)) \ra.\]
\end{Lem}

Since the coefficient $\la a^{(0)}_j(\epsilon), w_*^1 (a^{(1)}_i) \ra$
is given by the counting of monopoles
in $\M^{W_1} (a^{(1)}_i, a^{(0)}_j(\epsilon) )$ with the orientation,
and $\la a^{(0)}_j, w_*^0 (a^{(1)}_i(\epsilon)) \ra$ is given by the counting
of monopoles
in $\M^{W_0} (a^{(1)}_i(\epsilon), a^{(0)}_j )$ with the orientation,
the gluing Theorem \ref{final:glue}, for the admissible geometric limits
on $V\times \R$ and the corresponding holomorphic
triangles, yields the following Lemma (cf. Theorem 6.9 \cite{MW3}).

\begin{Lem}
\label{W_1:W_0} 
For small enough $\epsilon$ and large $r\ge r_0$, there 
is an orientation reversing  diffeomorphism  
\[   
\M^{W_1(r)}\bigl(a_i^{(1)}, a_j^{(0)}(\epsilon) \bigr) 
\cong \M^{W_0(r)} \bigl(a_i^{(1)}(\epsilon), a_j^{(0)}\bigr).
 \] 
Hence, we have $ w_*^0\circ w_*^1 =0$.
\end{Lem}

The Lemma corresponds to the fact that, in the two cases, the same
triangles are counted with the reverse orientation.

With these Lemmata at hand, the arguments in Section 6.5 \cite{MW3}
yield the following exact triangle for any $p \in \{0, \cdots, n-1\}$:
{\small
\ba
\diagram
HF^{SW}_{*} (Y_1, \s\otimes L_0) \rto^{w^1_{*}}
 & HF^{SW}_{*}(Y, \s\otimes L_0)
\dlto_{w^0_{*}} \\
 \bigoplus_{k\in \Z} HF^{SW}_{(*)}(Y_0, \s\otimes L_{nk+p})
\uto_{\Delta_{(*)}} 
& \\
\enddiagram
\label{triangle}
\na }

After a possible shift of relative index, we obtain the general
exact triangle for any $m, p \in \{0, \cdots, n-1\}$:
{\small
\[
\diagram
HF^{SW}_{*} (Y_1, \s\otimes L_m) \rto^{w^1_{*}}
 & HF^{SW}_{*}(Y, \s\otimes L_m)
\dlto_{w^0_{*}} \\
 \bigoplus_{k\in \Z} HF^{SW}_{(*)}(Y_0, \s\otimes L_{nk+p})
\uto_{\Delta_{(*)}} 
& \\
\enddiagram
\] }
This completes the proof of our main result, Theorem \ref{main:theorem},
except for the claim that  the connecting homomorphisms
 $\Delta_{(*)}$ are induced by the surgery cobordism
connecting $Y_0$ to $Y_1$. We shall discuss this statement in the next
subsection.

\subsection{The surgery triangle}

To give a precise description of the
connecting homomorphism $\Delta_{(*)}$, 
we need to study the discrepancy between 
the boundary operator
$\partial_Y$ of the Floer complex $C_*(Y,\s\otimes L_0, \mu)$
and the operator
$\partial_{Y_1}\oplus \bigoplus_{k\in \Z}
\partial_{Y_0,k}$ on 
$$ C_*(Y_1, \s\otimes L_0) \oplus 
\bigoplus_{k\in \Z} C_{(*)}(Y_0,\s\otimes L_{nk+p}). $$
Let us identify again the points of $\M_{Y,\mu}^*(\s\otimes L_0)$
as in (\ref{partition}) with 
$$ {\cal M}_{Y,\mu}^* (\s\otimes L_0)
= \{ a^{(1)}_i(\epsilon) \}_{i=1,\ldots r} \cup \{
a^{(0)}_j(\epsilon) \}_{j=r+1,\ldots, s }.$$
Then we have the following Lemma which
gives the connecting homomorphism $\Delta_{(*)}$.

\begin{Lem}
Suppose given a cycle in $\sum_i x_i a^{(0)}_i$ in
$C_{(*)}(Y_0,\s\otimes L_{nk+p})$. 
The image of $\sum_i x_i a^{(0)}_i$ under the
connecting homomorphism $\Delta$ is given by
$$ \Delta(\sum_i x_i a^{(0)}_i)= \sum_{i, j} x_i
\#\bigl(\M_{Y\times \R} (a_i^{(0)}(\epsilon), a_j^{(1)}(\epsilon))\bigr) 
a_j^{(1)}.$$
\end{Lem}
\begin{proof} (Lemma 7.1, Lemma 7.2 in \cite{MW3}) 
Using the diagram chasing,  Lemma \ref{inj:surj} and Lemma
\ref{W_1:W_0}, we see that we have
\[
 \begin{array}{lll}
&& \partial_{Y,\mu}(\sum_i x_i
a^{(0)}_i(\epsilon)) \\[2mm]
&=& \sum_{i,j} x_i 
\#\bigl(\M_{Y\times\R}(a_i^{(0)}(\epsilon),a_j^{(1)}(\epsilon))\bigr)
a_j^{(1)}(\epsilon) \\[2mm] 
&&- \sum_{i,j,k} x_i 
\#\bigl(\M_{Y\times\R}(a_i^{(0)}(\epsilon),a_j^{(1)}(\epsilon))\bigr)
\#\bigl(\M^{W_0}(a_j^{(1)}(\epsilon), a_k^{(0)}))\bigr)
a_k^{(0)}(\epsilon)\\[2mm]
&=& \sum_{i,j} x_i 
\#\bigl(\M_{Y\times\R}(a_i^{(0)}(\epsilon),a_j^{(1)}(\epsilon))\bigr)
a_j^{(1)}(\epsilon) \\[2mm] 
&&+ \sum_{i,j,k} x_i
\#\bigl(\M_{Y\times\R}(a_i^{(0)}(\epsilon),a_j^{(1)}(\epsilon))\bigr)
\#\bigl(\M^{W_1}(a_j^{(1)}, a_k^{(0)})(\epsilon))\bigr)
a_k^{(0)}(\epsilon).
\end{array} 
\]
By comparing this expression with
$$ 
\begin{array}{lll} 
&&w_*^1( \sum_i x_i \#\bigl( \M_{Y\times \R} (a_i^{(0)}(\epsilon),
a_j^{(1)}(\epsilon))\bigr) a_j^{(1)})  \\[2mm]
&=&\sum_j x_j \#\bigl( \M_{Y\times \R} (a_j^{(0)}(\epsilon),
a_i^{(1)}(\epsilon))\bigr) a_i^{(1)}(\epsilon)  \\[2mm]
&& +\sum_{i, j, k} x_i \#\bigl(\M_{Y\times \R}(a_i^{(0)}(\epsilon),
a_j^{(1)}(\epsilon))\bigr) \#\bigl(\M^{W_1}(a_j^{(1)}, a_k^{(0)}(\epsilon))
a_k^{(0)}(\epsilon)\bigr), 
\end{array} $$
we obtain 
$$ \Delta(\sum_i x_i a^{(0)}_i)= \sum_{i, j} x_i
\#\bigl(\M_{Y\times \R} (a_i^{(0)}(\epsilon), a_j^{(1)}(\epsilon)\bigr)
a_j^{(1)}.$$
This completes the proof. 
\end{proof}

In next proposition, we show that  the connecting homomorphism 
in the exact sequence can also be described as a map $w_{(*)}$,
induced by a surgery cobordism $\bar W_2$ connecting
$Y_0$ to $Y_1$, which satisfies $W_1 \#_Y W_0 = W_2 \# \bar{ \C\P^2}$.
The resulting diagram 
$$ C_*(Y_1, \s\otimes L_0) \stackrel{w_*^1}{\to} C_*(Y,\s\otimes L_0,\mu) 
\stackrel{w^0_*}{\to} \oplus_k 
C_{(*)}(Y_0,\s\otimes L_{nk+p}) \stackrel{w_{(*)}}{\to} C_*(Y_1,
\s\otimes L_0)[-1] $$
is therefore a distinguished triangle, the surgery triangle.  

\begin{Pro} (Proposition 7.3 \cite{MW3}) 
The connecting homomorphism $\Delta_{(*)}$
 in the exact triangle is given by 
the following expression,
$$ \Delta_{(*)}(\sum_i x_i a^{(0)}_i)= 
w_{(*)}(\sum_i x_i a^{(0)}_i),
$$
for any cycle $\sum_i x_i a^{(0)}_i$ in
$\oplus_k  C_{(*)}(Y_0,\s\otimes L_{nk+p})$ for any
fixed $p \in \{0, \cdots, n-1\},$ where
$$ w_{(*)} : \oplus_k C_{(*)}(Y_0,\s\otimes L_{nk+p})
\to C_*(Y_1, \s\otimes L_0)[-1]$$
is  the homomorphism defined by counting solutions in the
zero-dimensional components of the moduli spaces
$$ \M^{\bar W_2}(a_j^{(0)},a_i^{(1)}), $$
over the cobordism $\bar W_2$.
\end{Pro}
\begin{proof}
The argument is exactly the same as in the proof of Proposition
 7.3 in  \cite{MW3}, therefore we shall omit the proof here.
\end{proof}

\section{Seiberg--Witten and Casson--Walker invariant}

In this section, we derive the relation between the
topologically invariant version of the Seiberg-Witten invariant and
the Casson-Walker invariant for rational homology 3-spheres. Together
with the equivalence between the Casson-Walker invariant
and the theta invariant introduced by Ozsv\'ath and Szab\'o in
\cite{OS}, our result proves their  
conjecture relating the Seiberg-Witten invariant and their
theta invariant.

Let Y be a rational homology 3-sphere with a smoothly embedded knot $K$
representing a torsion element of order $n$ in $H_1(Y, \Z)$,
\[
\displaystyle{
\frac {|H_1(Y, \Z)|}{|\hbox{Torsion}( H_1(Y-\nu(K),  \Z))|}}
=n,
\]
and endowed with the canonical framing  $(m, l)$ in
a fixed identification:
$\nu (K)\cong D^2 \times S^1$. Let $p$ and $q$  be relatively prime
integers. The Dehn surgery with coefficient $p/q \in \Q \cup\{\infty\}$
on $K$ gives rise to another closed manifold $Y_{p/q}$.   

Denote by $\spinc (V)$ the set of equivalence classes of $\spinc$
structures on $V= Y\backslash \nu (K)$ with trivial restriction
to the boundary $T^2$. Then, for any $Y_{p/q}$, there
is a surjective map:
\[
\iota_{Y_{p/q}}: \qquad \spinc (Y_{p/q}) \to \spinc (V),
\]
where, for any $\s\in \spinc (Y_{p/q})$, $\iota_{Y_{p/q}} (\s)$
is given by the restriction to $V\subset Y_{p/q}$. The fiber of
$\iota_{Y_{p/q}}$ is given by a cyclic group
generated by the Poincar\'e dual of the core
of $Y_{p/q}\backslash V$.
Formally, for $\iota_{Y}(\s)$ with $\s \in \spinc (Y)$, we identify the fiber
of $\iota_{Y}$  with
the following set of $\spinc$ structures
\[
\iota_{Y}^{-1}(\iota_{Y}(\s)=\bigcup_{m=0, \cdots, n-1}
\{\s\otimes L_m| c_1(L_m) = mPD([K]) \in H^2(Y, \Z) \}.
\]
Similarly the fiber of $\iota_{Y_{p/q}}$ is given  by
\[
\bigcup_{m=0, \cdots, np-1} \{\s\otimes L_m|  c_1(L_m) 
= mPD([K]) \in H^2(Y_{p/q}, \Z) \},
\]
and the fiber of $\iota_{Y_{0}}$ is given  by
\[
\bigcup_{ m\in \Z}\{\s\otimes L_m| c_1(L_m) 
= mPD([K]) \in H^2(Y_0, \Z) \}.
\]
Here we use the same notation $\s$ on $Y_{p/q}$ ($Y$, or $Y_0$)
as the corresponding $\spinc$ structure obtained by gluing
$\s \in \spinc (V)$ with the trivial $\spinc$ structure on $\nu(K)$
by the trivial gauge transformation on $T^2$. We hope this
notation will not cause any confusion. 

Assume that $V$ and $\nu(K)$
are equipped with a metric with a cylindrical end modeled on $T^2$.
Let $\s$ be a $\spinc$ structure on $Y$. By the
result of \cite{CMW} on the moduli space of finite energy monopoles on
$(V, \iota_Y (\s))$, 
we know that
the irreducible part, denoted $\M^*_V(\s)$, is a smooth,
oriented 1-dimensional manifold. The 
asymptotic values along the cylindrical end define a
continuous map:
\[
\partial_\infty: \qquad \M^*_V(\s) \rightarrow  \chi_0 (T^2, V),
\]
where $\chi_0 (T^2, V)$ is a $\Z\times \Z_n$-covering of
the character torus $\chi (T^2)$. Sometime it is convenient
to compose the above asymptotic value map with this
covering map and define a boundary value map:
\ba
\partial_\infty: \qquad \M^*_V(\s) \rightarrow \chi (T^2).
\label{partial:infty}
\na

Notice that the reducible part $\chi (V)$ of the moduli space on
$(V, \iota_Y (\s))$
is an embedded circle $\chi (V) \subset \chi_0 (T^2, V)$ under the
asymptotic value map. This becomes a circle of 
multiplicity $n$ in $\chi (T^2)$. There is a ``bad point'' 
in $\chi (T^2)$, given by the flat connections such that the
corresponding twisted Dirac operator has a non-trivial kernel.
We can endow $\chi (T^2)$ with a coordinate system $(u, v)$ defined
by the holonomy around the longitude $l$ and the meridian $m$,
respectively, so that the bad point corresponds to $(u, v) =(1,
1)$. The reducible circle $\chi (V)$, with the holonomy around the
longitude $l$ of order $n$, is given by $u = u(\s)$, with
$u(\s) \in \{0, 2/n, \cdots, 2(n-1)/n\}$. After a suitable perturbation,
and a corresponding shift of coordinates, as discussed in \cite{CMW},
we can assume that the bad point does not lie on any of these $n$
possible circles $u = u(\s)$ of reducibles $\chi (V)$. 

 From the result in \cite{CMW}, we know
that, under the map $\partial_\infty$ in (\ref{partial:infty}), the
boundary points $\partial ( \M^*_V(\s))$ are
either mapped to the bad point in $\chi(T^2)$ or mapped
to the reducible circle $u= u(\s)$ on $\chi(T^2)$. 

Let $\chi (\nu(K) \subset Y_{p/q})$ be the reducible circle
on $\nu(K) \subset Y_{p/q}$, which maps to a closed curve on $\chi(T^2)$
with slope $p/q$ in the $(u, v)$-coordinates: parallel to $pv =qu$. 
Looking at the induced Spin structure on $T^2 \subset Y_{p/q}$,
we know that the curve $\chi (\nu(K) \subset Y_{p/q})$ goes through
$(0, 1)$ if $q$ is odd or goes through
$(0, 0)$ if $q$ is even, cf.\cite{CMW}. Again, after a suitable
perturbation as in (\ref{eta:pert}), and the corresponding shift of 
coordinates, we can assume that this $p/q$-curve is away from the bad 
point on $\chi(T^2)$ and does not meet $u= u(\s)$ 
along the coordinate line $v=0$. Then we know that
$u= u(\s)$ intersects $\chi (\nu(K) \subset Y_{p/q})$ 
inside $\chi(T^2)$ at $p$ points, which
are denoted by $\theta_1, \cdots, \theta_p$, ordered according the
orientation of $ u= u(\s) \subset \chi (T^2)$. They can be
lifted to $pn$ points in $\chi_0(T^2, V)$. We denote these points by
$\theta_1^{(k)},  \cdots,  
\theta_p^{(k)}$, ($k=0, 1, \cdots, n-1$)
according to the order. Denote by $\theta_0$ the intersection point of 
$u= u(\s)$ with $v=0$ in $\chi (T^2)$. This
can be lifted to $n$-points $\theta_0^{(0)}, 
\theta_0^{(1)}, \cdots, \theta_0^{(n-1)}$
on $\chi_0(T^2, V)$.  
Moreover, we can assume that the map $\partial_\infty$ in
(\ref{partial:infty}) 
is transverse to the curves $u= u(\s)$, $v=0$ and
$\chi (\nu(K) \subset Y_{p/q})$, by a suitable perturbation of the
Seiberg-Witten equations on $V$ as in \cite{CMW}. We can also assume
that the image   
$\partial_\infty\M^*_V(\s)$ does not meet the points $\theta_0,
\theta_1, \cdots \theta_p$ in $\chi (T^2)$, again
by suitable perturbation, as discussed in \cite{CMW}.

Let $I$ be any open interval in
\[
 \chi (V) = \{ u = u(\s)\} \subset \chi_0(T^2, V).
\]
We denote by $SF_\C (\dirac^V_I)$
the complex spectral flow of Dirac operator
on $V$ twisted with the path of reducible connections $I$ on $V$. From  
the analysis in \cite{CMW} and \cite{MW1},
we know that 
\ba
\# \bigl( \partial_\infty|_{\partial \M^*_V(\s)} \bigr)^{-1}
(I) = SF_\C (\dirac^V_I).
\label{boundary:on:I}
\na
   For convenience, we define
\ba
SF_\C (\dirac^V_{[\theta_i, \theta_j]}) = \sum_{k=0}^{n-1}
SF_\C (\dirac^V_{[\theta_i^{(k)}, \theta_j^{(k)}]}).
\label{spectralflow:sum}
\na

With this notation understood, we can state the following
proposition relating the Seiberg-Witten invariants on $Y_{p/q}$,
$Y$ and $Y_0$. 

\begin{Pro}\label{p/q-formulae:SW} 
Consider generic compatible small perturbations of the Seiberg-Witten
equations 
on $Y_{p/q}$, $Y$ and $Y_0$, such that the map $\partial_\infty$ as in
(\ref{partial:infty}) 
is transverse to the curves $u= u(\s)$, $v=0$ and
$\chi (\nu(K) \subset Y_{p/q})$ and misses the points
$\theta_0, \theta_1, \cdots, \theta_p$ in $\chi (T^2)$. Then
we have the following relation:
\[\begin{array}{lll}
&&\sum_{k=0}^{pn-1} SW_{Y_{p/q}} (\s\otimes L_k, g_{Y_{p/q}})\\[2mm]
&=& p\sum_{k=0}^{n-1} SW_{Y} (\s\otimes L_k, g_{Y}) +
q \sum_{k\in \Z} SW_{Y_0} (\s\otimes L_k)\\[2mm]
&& + \sum_{i=1}^p SF_\C (\dirac^V_{[\theta_0, \theta_i]}).
\end{array}
\]
\end{Pro}
\begin{proof} By the gluing theorem for 3-dimensional monopoles as in
\cite{CMW}, we have 
\[
\bigcup_{k=0}^{pn-1} \M^*_{Y_{p/q}} (\s\otimes L_k) 
= \M^*_V (\s) \times_{\chi (T^2)} \chi (\nu(K) \subset Y_{p/q}),
\]
where $\chi (\nu(K) \subset Y_{p/q})$ is the $p/q$-curve on 
$\chi (T^2)$. Thus, we obtain 
\ba
\sum_{k=0}^{pn-1} SW_{Y_{p/q}} (\s\otimes L_k, g_{Y_{p/q}})
= \# \bigl( \M^*_V (\s) \times_{\chi (T^2)} \chi (\nu(K) \subset Y_{p/q})
\bigr).
\label{p:q}
\na
Notice that the set 
$
\{\theta_1^{(k)}, \cdots, \theta_p^{(k)}: k = 0, 1, \cdots, n-1 \}
$
consists of  the unique reducible monopole for each 
$(Y_{p/q}, \s\otimes L_k)$.

Similarly, we have
\ba
\sum_{k=0}^{n-1} SW_{Y} (\s\otimes L_k, g_{Y})
=  \# \bigl( \M^*_V (\s) \times_{\chi (T^2)} \{ v=0\}\bigr).
\label{Y}
\na
Here the reducible set consists of
$
\{ \theta_0^{(0)}, \theta_0^{(1)}, \cdots, \theta_0^{(n-1)} \}.
$

In order to avoid the circle of reducibles on $(Y_0, \s \otimes L_0)$,
we need to introduce a small perturbation such that
$\chi (\nu(K)\subset Y_0)$ on $\chi (T^2)$ is a small parallel
shifting of $u = u(\s)$ such that the bad point is not
contained in the narrow strip bounded by these
two parallel curves. We denote this small shift
of $u = u(\s)$ by $u = u(\s) +\eta$, where $\eta$ is a sufficiently small
positive number. This can be achieved by a perturbation of the
equations as in \cite{CMW}. Then we have 
\ba
\sum_{k\in \Z}  SW_{Y_0} (\s\otimes L_k) 
= \# \bigl( \M^*_V (\s) \times_{\chi (T^2)} \{ u=  u(\s) +\eta\}\bigr).
\label{Y_0}
\na

In order to compare the three countings in (\ref{p:q}) -- (\ref{Y_0}), we
need to choose an oriented 2-chain $ C$ in $\chi (T^2)$ whose
boundary 1-chain is given by 
\[\begin{array}{lll}
&&\chi (\nu(K)\subset Y_{p/q}) - p \chi (\nu(K)\subset Y) -q
\chi (\nu(K)\subset Y_0) \\[2mm]
&=& \chi (\nu(K)\subset Y_{p/q}) - p \{ v= 0\} 
- q\{ u= u(\s) +\eta\},
\end{array}
\]
and such that $C$ does not contain the bad point in $\chi (T^2)$. Then, 
counting the boundary points of $\partial_\infty^{-1}(C)$, as a
0-chain, we obtain
\ba
\begin{array}{lll}
&&\# \bigl( \partial_\infty^{-1} (\chi (\nu(K)\subset Y_{p/q}) \bigr)
\\[2mm]
&=& p \# \bigl( \partial_\infty^{-1} (\{ v= 0\}) \bigr)
 + q  \# \bigl( \partial_\infty^{-1}(\{ u= u(\s) +\eta\}) \bigr)\\[2mm]
&& + \# \bigl( \partial_\infty\mid_{\partial (\M^*_V(\s)) }\bigr)^{-1} (C).
\label{count:C}
\end{array}
\na

As $C$ does not contain the bad points, we know
that the possible points of $\partial_\infty (\partial (\M^*_V(\s)) )
\cap C$ all lie on the curve $u = u(\s)$, away from
the points $\theta_0, \theta_1, \cdots, \theta_p$. It is easy to see
that $C$ covers the intervals of $u = u(\s)$ between two consecutive
points $\theta_i$ with different multiplicities:
the multiplicities are 
$ p, p-1, \cdots, 1, 0$,
for the intervals 
\[
[\theta_0, \theta_1], [\theta_1, \theta_2],
\cdots, [\theta_{p-1}, \theta_p], [\theta_p, \theta_0],
\]
respectively. 
By the identity (\ref{boundary:on:I}) and the definition 
(\ref{spectralflow:sum}), we know that 
\ba
\# \bigl( \partial_\infty|_{\partial (\M^*_V(\s))} \bigr)^{-1}
  (C) = \sum_{i=1}^p SF_\C (\dirac ^V_{[\theta_0, \theta_i]}).
\label{count:u(s)}
\na
Combining all the identities in (\ref{p:q}), (\ref{Y}),
(\ref{Y_0}),  (\ref{count:C}) and (\ref{count:u(s)}),
we obtain the proof of the proposition.
\end{proof}

The Seiberg-Witten invariant for any rational homology 3-sphere
depends on metric and perturbation (cf.\cite{MW1}). We now consider 
the correction term (\ref{correction:term}) as defined in the
introduction. We have the following proposition relating the
correction terms for $Y_{p/q}$ and $Y$.

\begin{Pro}\label{p/q-formulae:correction}
\begin{enumerate}
\item For any rational homology 3-sphere $Y$ with a $\spinc$ structure
$\s$ and a Riemannian metric $g_Y$, 
\[
\hat{SW}_Y (\s) = SW_Y (\s, g_Y) - \xi (\s, g_Y)
\]
is a well-defined topological invariant.
\item For any relatively prime integers $p$ and $q$, a positive integer $n$,
and $u\in \{0, 2/n,\ldots, 2(n-1)/n \}$, we have that
 \[
\sum_{k=0}^{pn-1} \xi_{Y_{p/q}} (\s\otimes L_k, g_{Y_{p/q}})
- p \sum_{k=0}^{n-1} \xi_{Y} (\s\otimes L_k, g_{Y}) -
\sum_{i=1}^p SF_\C (\dirac ^V_{[\theta_0, \theta_i]}) \]
 is independent of the 
manifold $Y$ and depends only on $p, q, n$, and
$u(\s) \in \{0, 2/n, \ldots, 2(n-1)/n\}.$
\end{enumerate}
\end{Pro} 
\begin{proof} Claim (1) follows from the wall-crossing formulae
in \cite{MW1} and the Atiyah-Patodi-Singer index theorem. 
The proof of claim (2) is analogous to the proof of Proposition
7.9 in \cite{OS}. We adapt their arguments to our situation.
We write the standard surgery cobordism between
$S^3$ and the Lens space $L(p, q)$ as
\[
W(S^3, L(p, q)) = \bigl( [0, 1] \times S^1\times D^2 \bigr) 
\cup_{ [0, 1] \times S^1 \times S^1} X_{p/q},
\]
Then the surgery cobordism  between $Y$ and $Y_{p/q}$
can be identified as 
\[
W_{p/q} =  \bigl( [0, 1] \times V\bigr) 
\cup_{ [0, 1] \times S^1 \times S^1} X_{p/q}.
\]
We fix a metric on $W_{p/q}$ which respects the product structure
$[0, 1] \times V$ and $[0, 1] \times S^1 \times S^1$,
and agrees with $g_Y$ and $g_{Y_{p/q}}$ on the boundaries
$Y$ and $Y_{p/q}$, respectively.

For a $\spinc$ structure  $\s\otimes L_i^{(m)}$
in $\{\s\otimes L_k: k=0, \cdots, pn-1\}$ on
$Y_{p/q}$, whose reducible monopole corresponds
to $\theta_i^{(m)}$ (with $i\in \{1, \cdots, p\}$ and
$m\in \{ 0, \cdots n-1\}$), we consider the $\spinc$ structure
$\s\otimes L_m$ on $Y$ whose  reducible monopole is
$\theta_0^{(m)}$. Then we claim that 
\ba\label{xi:m:i}
\xi_{Y_{p/q}} (\s\otimes L_i^{(m)}, g_{Y_{p/q}}) 
- \xi_{Y} (\s\otimes L_m, g_Y) -SF_\C
(\dirac^V_{[\theta_0^{(m)}, \theta_i^{(m)}]})\na
is independent of $Y$ and depends only on $p, q, n$ and on
$u(\s) \in \{0, 2/n, \cdots, 2(n-1)/n\}.$

To prove this claim, we choose a $\spinc$ structure $\tilde \s$
on $W_{p/q}$ whose restriction to $Y$ and $Y_{p/q}$ is given by
$\s\otimes L_m$ and $\s\otimes L_i^{(m)}$, respectively, and such that
$c_1(\tilde\s)^2 =1$. On $(W_{p/q}, \tilde \s)$, we choose
a connection $A$, whose restriction to
$V\times [0, 1]$ is the path of reducibles connecting
$\theta_0^{(m)}$ to $\theta_i^{(m)}$ along the curve
$\chi (V) \subset \chi_0(T^2, V)$. Then we have
\ba\begin{array}{lll}
&& \xi_{Y_{p/q}} (\s\otimes L_i^{(m)}, g_{Y_{p/q}}) 
- \xi_{Y} (\s\otimes L_m, g_Y) \\[2mm]
&=& Ind_\C (\Dirac_A^{W_{p/q}}) - \displaystyle {
\bigl( \frac {c_1(\tilde \s)^2 -\sigma (W_{p/q}) }{8}\bigr)}\\[2mm]
&=& Ind_\C (\Dirac_A^{W_{p/q}})\\[2mm]
&=& Ind_\C (\Dirac_A^{[0,1]\times V}) + Ind_\C (\Dirac_A^{X_{p/q}})
\end{array}\label{splitting}
\na
where the third equality follows from the
splitting principle for the index, as the
Dirac operator has no kernel on the various boundaries and
corners \cite{CLM}. Notice that we have
\[
Ind_\C (\Dirac_A^{[0,1]\times V}) =
SF_\C (\dirac^V_{[\theta_0^{(m)}, \theta_i^{(m)}]}),
 \]
and the connection $A|_{X_{p/q}}$ extends to connection $A_0$
on $W(S^3, L(p, q))$ by a flat connection, whose index on
$[0,1] \times S^1 \times D^2$ satisfies
\[
Ind_\C (\Dirac_{A_0}^{[0,1] \times S^1 \times D^2} )=0.
\]
In fact, we can choose the metric on $W(S^3, L(p, q))$ with
a positive scalar curvature metric on $[0,1] \times S^1 \times D^2$.
Therefore, we have
\[
Ind_\C (\Dirac_A^{X_{p/q}}) = 
Ind_\C (\Dirac_{A_0}^{W(S^3, L(p, q))}),
\]
which depends only on $p, q, n$ and $u(\s)$, and so does the quantity
\ba
\label{splitting2}
 \xi_{Y_{p/q}} (\s\otimes L_i^{(m)}, g_{Y_{p/q}})
- \xi_{Y} (\s\otimes L_m, g_Y)-
SF_\C (\dirac^V_{[\theta_0^{(m)}, \theta_i^{(m)}]}).
\na 
When summing the identity (\ref{splitting2})
over $i\in \{1, \cdots, p\}$ and $m\in \{0, \cdots, n-1\}$, 
notice that the term $\xi_{Y} (\s\otimes L_m, g_Y)$ is independent of 
$i\in \{1, \cdots, p\}$, 
hence we obtain the proof of the claim (2) by using the definition 
(\ref{spectralflow:sum}).  
\end{proof}

With these two propositions in place, we now have the following surgery
formula for the modified version of the Seiberg-Witten invariant.

\begin{The}\label{p/q-formulae:hatSW}
Given any two relatively prime integers $p$ and $q$, a positive
integer $n$ and $u\in \{0, 2/n, 2(n-1)/n\}$, there is a rational valued
function $s(p, q, n, u)$, depending only on $p, q, n$ and $u$,
satisfying the
following property. Let $Y$ be a rational homology 3-sphere with a smoothly
embedded knot and a canonical framing $(m, l)$
such that $\nu (K)\cong D^2\times S^1$. Assume that
$K$ represents a torsion element of order $n$ in $H_1(Y, \Z)$. Let
$\s$ be a $\spinc$ structure on $Y$. Then we have
\[
\begin{array}{lll}
&&\sum_{k=0}^{pn-1} \hat{SW}_{Y_{p/q}} (\s\otimes L_k)\\[2mm]
&=& p\sum_{k=0}^{n-1} \hat{SW}_{Y} (\s\otimes L_k) +
q \sum_{k\in \Z} SW_{Y_0} (\s\otimes L_k)\\[2mm]
&& + s(p, q, n, u). \end{array}
\]
\end{The}
\begin{proof} Following from Proposition \ref{p/q-formulae:correction}, 
we know that
\ba
\sum_{k=0}^{pn-1} \xi_{Y_{p/q}} (\s\otimes L_k, g_{Y_{p/q}})
-p \sum_{k=0}^{n-1} \xi_{Y} (\s\otimes L_k, g_{Y}) -
\sum_{i=1}^p SF_\C (\dirac ^V_{[\theta_0, \theta_i]})
\label{s(p,q,n,u)}
\na
 depends only on $p, q, n$ and
$u=u(\s) \in \{0, 2/n, \cdots, 2(n-1)/n\}$. We denote this term by
$s(p, q, n, u)$. By subtracting (\ref{s(p,q,n,u)}) from the surgery
formula for the Seiberg-Witten invariants in Proposition  
\ref{p/q-formulae:SW}, we obtain the proof of this theorem.
\end{proof}

Now we can establish the equivalence between
the modified version of the Seiberg-Witten invariant $\hat{SW}$
and the Casson-Walker invariant for rational homology 3-spheres.

\begin{The}
For any rational homology 3-sphere $Y$, we have
\[
\sum_{\s\in \spinc (Y)} \hat {SW}_Y (\s) 
= \displaystyle{\frac 12} |H_1(Y, \Z)| \lambda (Y)
\]
where $\lambda (Y)$ is the Casson-Walker invariant.
\end{The}
\begin{proof}
We first derive the surgery formula for the invariant
$\sum_{\s\in \spinc (Y)} \hat {SW}_Y (\s)$ from 
Theorem \ref{p/q-formulae:hatSW} and the Seiberg-Witten invariant for
$Y_0$ (a rational homology $S^1\times S^2$, i.e., $b_1(Y_0)=1$) (see
\cite{MT} \cite{Don}):
\ba\label{surgery:formuale:SW}
\begin{array}{lll}
&&\sum_{\s\in \spinc (Y_{p/q})} 
\hat{SW}_{Y_{p/q}} (\s)\\[2mm]
&=& p\sum_{\s\in \spinc (Y)} \hat{SW}_{Y} (\s) +
q \sum_{j=0}^\infty a_j j^2 + |H_1(Y, \Z)| s(p, q, n) \end{array}
\na
where $s(p, q, n)= \sum_u s(p, q, n, u)/n$ and $a_j$ is the coefficient
of the symmetrized Alexander polynomial of $Y_0$,
\[
A(t) = |\hbox{Torsion} (H_1(Y_0, \Z))| + \sum_{j=1}^\infty a_j (t^j + t^{-j})
\]
normalized such that 
\[
A(1) =  |\hbox{Torsion} (H_1(Y_0, \Z))|.
\]
Set $\bar \lambda (Y) = \displaystyle{\frac 12}|H_1(Y, \Z)| \lambda (Y)$
as the normalized Casson-Walker invariant. Then
the surgery formula in \cite{Walker}
for $\bar \lambda (Y)$ can be expressed as (cf. \cite{OS}):
\ba\label{surgery:formuale:CW}
\begin{array}{lll}
\bar \lambda (Y_{p/q}) &=& p \bar \lambda (Y) + q  \sum_{j=0}^\infty a_j j^2 
\\[2mm]
&& + |H_1(Y, \Z)| \bigl( \displaystyle{
\frac{q(n^2-1)}{12n^2} -\frac{p s(p, q)}{2}}\bigr).
\end{array}\na
Here $s(p, q)$ is the Dedekind sum of relatively prime integers $p$
and $q$ (cf. \cite{Walker}). 
Comparing (\ref{surgery:formuale:SW}) and (\ref{surgery:formuale:CW}),
we only need to show that 
\ba
s(p, q, n) =  \displaystyle{
\frac{q(n^2-1)}{12n^2} -\frac{p s(p, q)}{2}}.
\label{s(p,q,n)}
\na

Since $s(p,q,n)$ is independent of the manifold $Y$, we can choose
some examples that can be computed explicitly, and use them to
identify the coefficient $s(p,q,n)$.  
The Lens space $L(p, q)$ can be obtained by a $p/q$-surgery
on an unknot in $S^3$. The calculation of Nicolaescu \cite{Nic}
for $L(p, q)$ gives us that 
\[
\sum_{\s\in \spinc (L(p, q))} \hat{SW}_{L(p, q)}(\s) = 
-\displaystyle{
\frac{ps(p, q)}{2}}.
\]
This implies that (\ref{s(p,q,n)}) holds for $n=1$. Now
we can prove (\ref{s(p,q,n)}) by induction on $n$. 
This is exactly the same argument as in the proof of Theorem 7.5 in
\cite{OS} on the equivalence of their theta invariant and the
Casson-Walker invariant. The example is the Seifert manifold
$M(n, 1; -n, 1; q, -p)$, obtained by $p/q$ surgery on a knot of order $n$ in 
$L(n, 1) \# \overline{L(n, 1)}$. By Kirby calculus it is possible to show
that $M(n, 1; -n, 1; q, -p)$ can be obtained as
$(-n)$-surgery on a knot in the Lens space $L(pn-q, q)$, and
can be obtained as a sequence of surgeries on knots of order less than $n$,
see the proof of Theorem 7.5 in \cite{OS} for details. 
\end{proof}

\noindent {\bf Matilde Marcolli}\par 
\noindent Max--Planck--Institut f\"ur Mathematik, Vivatsgasse 7,
D-53111 Bonn, Germany. \par 
\noindent marcolli\@@mpim-bonn.mpg.de

\vskip .2in

\noindent {\bf Bai-Ling Wang}\par
\noindent  Department of Pure
Mathematics, University of Adelaide, Adelaide SA 5005, Australia. \par
\noindent bwang\@@maths.adelaide.edu.au

\end{document}